\documentclass[10pt]{amsart}
\usepackage{amsmath}
\usepackage{amssymb}
\usepackage{amsfonts}
\usepackage{amsthm}
\usepackage{amsxtra}
\usepackage{rotating}
\usepackage{nicefrac}
\usepackage{float}
\usepackage[all,web,arc,poly,dvips]{xy}
\usepackage{epic,eepic}
\usepackage{epsfig}
\usepackage[dvips]{color}
\usepackage{array}
\usepackage{pifont}
\usepackage{multirow}

\theoremstyle{plain}
\newtheorem{theorem}{Theorem}[section]

\newtheorem{lemma}{Lemma}[section]
\newtheorem{proposition}{Proposition}[section]

\theoremstyle{definition}

\theoremstyle{remark}
\newtheorem{remark}{Remark}[section]

\newcommand{\mfS}{\mathfrak S}
\newcommand{\mfM}{\mathfrak M}
\newcommand{\CC}{\mathbb C}
\newcommand{\CP}{\mathbb C\mathbb P}
\newcommand{\NN}{\mathbb N}
\newcommand{\ZZ}{\mathbb Z}

\newcommand{\pII}{${\rm P}_{\rm II}\:$}
\newcommand{\PV}{${\rm P}_{\rm V}\:$}
\newcommand{\dPV}{${\rm dP}_{\rm V}\:$}
\newcommand{\PVI}{${\rm P}_{\rm VI}\:$}

\newcommand{\half}{
        {\lower0.00ex\hbox{\raise.6ex\hbox{\the\scriptfont0 1}
                           \kern-.5em\slash\kern-.1em\lower.45ex
                                     \hbox{\the\scriptfont0 2}}}}
\newcommand{\quarter}{
        {\lower0.00ex\hbox{\raise.6ex\hbox{\the\scriptfont0 1}
                           \kern-.5em\slash\kern-.1em\lower.45ex
                                     \hbox{\the\scriptfont0 4}}}}
\newcommand{\tquarter}{
        {\lower0.00ex\hbox{\raise.6ex\hbox{\the\scriptfont0 3}
                           \kern-.5em\slash\kern-.1em\lower.45ex
                                     \hbox{\the\scriptfont0 4}}}}
\newcommand{\eighth}{
        {\lower0.00ex\hbox{\raise.6ex\hbox{\the\scriptfont0 1}
                           \kern-.5em\slash\kern-.1em\lower.45ex
                                     \hbox{\the\scriptfont0 8}}}}
\newcommand{\othird}{
        {\lower0.00ex\hbox{\raise.6ex\hbox{\the\scriptfont0 1}
                           \kern-.5em\slash\kern-.1em\lower.45ex
                                     \hbox{\the\scriptfont0 3}}}}

\begin{document}

\title[Random Matrix Theory and the Sixth Painlev\'e Equation]
{Random Matrix Theory and the Sixth Painlev\'e Equation}

\author{P.J.~Forrester and N.S.~Witte}
\address{Department of Mathematics and Statistics,
University of Melbourne, Victoria 3010, Australia}
\email{\tt p.forrester@ms.unimelb.edu.au; n.witte@ms.unimelb.edu.au}

\vspace{5mm}
\begin{center}
{Dedicated to the centenary of the publication of the Painlev\'e VI equation\\
 in the Comptes Rendus de l'Academie des Sciences de Paris\\
 by Richard Fuchs in 1905.}
\end{center}
\vspace{5mm}

\begin{abstract}
A feature of certain ensembles of random matrices is 
that the corresponding measure is invariant under
conjugation by unitary matrices. Study of such ensembles
realised by matrices with Gaussian entries leads to
statistical quantities related to the eigenspectrum, such
as the distribution of the largest eigenvalue, which can
be expressed as multidimensional integrals or equivalently
as determinants. These distributions are well known to be
$\tau$-functions for Painlev\'e systems, allowing for the
former to be characterised as the solution of certain
nonlinear equations. We consider the random matrix ensembles
for which the nonlinear equation is the $\sigma$ form of
\PVI. Known results are reviewed, as is their
implication by way of series expansions for the distributions.
New results are given for the boundary conditions in the 
neighbourhood of the fixed singularities at $t=0,1,\infty$ of
$\sigma$\PVI displayed by a generalisation of the generating
function for the distributions. The structure of these expansions
is related to Jimbo's general expansions for the $\tau$-function
of $\sigma$\PVI in the neighbourhood of its fixed singularities,
and this theory is itself put in its context of the linear
isomonodromy problem relating to \PVI.
\end{abstract}

\maketitle

\section{Introduction}
\setcounter{equation}{0}

\subsection{The $\sigma$-form of Painlev\'e VI}
Given a large sequence of an eigenvalue spectrum, it is a simple matter to
rescale so that the mean spacing between consecutive eigenvalues is unity, then 
to empirically determine the distribution function for the spacing. When these 
eigenspectra are the highly excited states of heavy nuclei, it is a celebrated 
result that the distribution function is well approximated by the functional 
form
\begin{equation}
  p^{(\rm W)}_1(s) := \frac{\pi s}{2}e^{-\pi s^2/4} , 
\end{equation}
known as the Wigner surmise. It is in the exact computation of eigenvalues 
distributions for certain classical random matrix ensembles that Painlev\'e
transcendents make their appearance. These transcendents may relate to any of 
\pII to \PVI, depending on the random matrix ensemble under consideration,
or its scaled limits. As part of the centenary of the discovery of Painlev\'e
VI by Fuchs \cite{Fu_1905}, in this paper we will restrict attention to the
random matrix ensembles relating to \PVI.

Painlev\'e VI conventionally refers to the four parameter second order 
nonlinear differential equation
\begin{multline}
  y'' = \frac{1}{2}\left(\frac{1}{y}+\frac{1}{y-1}+\frac{1}{y-t} \right)(y')^2
        -\left(\frac{1}{t}+\frac{1}{t-1}+\frac{1}{y-t} \right)y' \\
        +\frac{y(y-1)(y-t)}{t^2(t-1)^2}
         \left(\alpha+\frac{\beta t}{y^2}+\frac{\gamma(t-1)}{(y-1)^2}
                     +\frac{\delta t(t-1)}{(y-t)^2} \right) ,
\label{PVI_y}
\end{multline} 
as obtained by Fuchs \cite{Fu_1905}. However in random matrix theory this equation
is never encountered directly. Rather what is encountered is the so-called 
Jimbo-Miwa-Okamoto $\sigma$-form of Painlev\'e VI,
\begin{multline}
\sigma'_{\rm VI}\left(t(t-1)\sigma''_{\rm VI}\right)^2
  + \left( \sigma'_{\rm VI}\left[2\sigma_{\rm VI}-(2t-1)\sigma'_{\rm VI}\right]
             +v_1v_2v_3v_4 \right)^2 = \prod^4_{k=1}(\sigma'_{\rm VI}+v^2_k) .
\label{PVI_sigmaF}
\end{multline}
This is written in a form which displays a $ D_4 $ root system symmetry in the 
parameters $ v_1,\ldots,v_4 $. When expanded out there is a common factor of 
$ \sigma'_{\rm VI} $, which when cancelled out shows (\ref{PVI_sigmaF}) to be 
a second order second degree nonlinear differential equation. The equations 
(\ref{PVI_y}) and (\ref{PVI_sigmaF}) are
related by the Hamiltonian formulation of \PVI, due originally to 
Malmquist in 1922 \cite{Ma_1922}. 

In the Hamiltonian approach to the Painlev\'e equations in general, one presents
a Hamiltonian $ H(q,p,t;\{v_k\}) $ where $ \{v_k\} $ are parameters, such that
after eliminating $ p $ in the Hamilton equations
\begin{equation}
   q' = \frac{\partial H}{\partial p}, \qquad
   p' =-\frac{\partial H}{\partial q} ,
\label{HamEq}
\end{equation}
where the dash denotes derivatives with respect to $ t $, the equation in $ q $
is the appropriate Painlev\'e equation. The Hamiltonians can be systematically 
derived from the isomonodromy deformation theory associated with the Painlev\'e
equations \cite{KO_1984}, \cite{Ok_1986a} (aspects of the isomonodromy deformation 
theory associated with \PVI is covered in Section \ref{sectionPVI_isoMD} below). 
From such considerations, the 
Hamiltonian relating to \PVI was given by Okamoto \cite{Ok_1987a} as
\begin{multline}
  t(t-1)H_{\rm VI} = q(q-1)(q-t)p^2 \\
  -\left[ (v_3+v_4)(q-1)(q-t)+(v_3-v_4)q(q-t)-(v_1+v_2)q(q-1) \right]p \\
  +(v_3-v_1)(v_3-v_2)(q-t)
\label{PVI_Ham}
\end{multline}
where the the parameters $ v_1,\ldots,v_4 $ are related to 
$ \alpha,\beta,\gamma,\delta $ in (\ref{PVI_y}) according to
\begin{equation}
  \alpha = \frac{1}{2}(v_1-v_2)^2, \quad
  \beta  =-\frac{1}{2}(v_3+v_4)^2, \quad
  \gamma = \frac{1}{2}(v_3-v_4)^2, \quad
  \delta = \frac{1}{2}(1-(1-v_1-v_2)^2),
\end{equation}
and $ q $ satisfies (\ref{PVI_y}). Note that $ H_{\rm VI} $ is quadratic in 
$ p $ and thus according to the first of the Hamilton equations (\ref{HamEq})
$ p $ can be written as a rational function of $ q $ and $ q' $, and thus in
fact $ H_{\rm VI} $ is a rational function in $ q $ and $ q' $. According to
the following result, this particular rational function, augmented by the 
addition of a linear function in $ t $, satisfies (\ref{PVI_sigmaF}) 
\cite{Ok_1987a}.

\begin{proposition}
Define the auxiliary Hamiltonian 
\begin{equation}
  h_{\rm VI}(t) = t(t-1)H_{\rm VI} 
  +e_2[-v_1,-v_2,v_3]t-\frac{1}{2}e_2[-v_1,-v_2,v_3,v_4] ,
\end{equation}
where 
\begin{equation}
  e_p[a_1,\ldots,a_s] := \sum_{1\leq j_1<\ldots<j_p\leq s}a_{j_1}a_{j_2}\cdots a_{j_p} .
\end{equation}
This auxiliary Hamiltonian satisfies the $\sigma$-form of Painlev\'e VI 
(\ref{PVI_sigmaF}).
\end{proposition}

\subsection{Historical Overview}
There are certain ensembles of $ N\times N $ random matrices with complex Gaussian
entries and invariance under conjugation by unitary matrices, 
which have their joint eigenvalue probability distribution function (p.d.f.)
of the form
\begin{equation}
   \frac{1}{C}\prod^{N}_{l=1}w_2(x_l)\prod_{1\leq j<k\leq N}(x_k-x_j)^2 ,
\label{Epdf}
\end{equation}
where the weight function $ w_2(x) $ is of one of the classical 
forms
\begin{equation}
   w_2(x) = \begin{cases}
             e^{-x^2}, & {\rm Gaussian} \\
             x^ae^{-x} \quad (x>0), & {\rm Laguerre} \\
             x^a(1-x)^b \quad (0<x<1), & {\rm Jacobi}
            \end{cases} .
\label{classicalW}
\end{equation}
For example, let $ X $ be a $ n\times N (n\geq N) $ rectangular matrix of complex
Gaussians $ {\rm N}[0,1/\sqrt(2)]+i{\rm N}[0,1/\sqrt(2)] $. Then the matrix 
$ X^{\dagger}X $ has eigenvalue p.d.f. (\ref{Epdf}) with Laguerre weight 
$  x^{n-N}e^{-x} \quad (x>0) $. The sixth Painlev\'e equation relates to 
(\ref{Epdf}) with Jacobi weight; in particular to the probability that there are
exactly $ n $ eigenvalues in the interval $ (t,1) $ of that ensemble. This
probability is in turn equal to the coefficient of $ (1-\xi)^n $ in the expansion
of 
\begin{multline}
  E^J_N(t;a,b;\xi) 
  := \frac{1}{C}\left(\int^1_0-\xi\int^1_t\right)dx_1 \cdots
                \left(\int^1_0-\xi\int^1_t\right)dx_N \\
     \times \prod^N_{l=1}x^a_l(1-x_l)^b \prod_{1\leq j<k \leq N}(x_k-x_j)^2 ,
\label{EJUE}
\end{multline}
where $ C $ denotes the normalisation.
In the case $N=1$ (\ref{EJUE}) is an integral form of a particular
${}_2F_1$ hypergeometric function, which was related to
the Painlev\'e VI equation by Okamoto \cite{Ok_1987a}.

Let $ \{p_j(x)\}_{j=0,1,\ldots} $ denote the set of monic polynomials of degree
$ j $ orthogonal with respect to the Jacobi weight
$ x^a(1-x)^b (0<x<1) $, which are given in terms of the Jacobi polynomials 
$ P^{(a,b)}_j $ by
\begin{equation}
   p_j(x) = (-1)^j j!\frac{\Gamma(a+b+j+1)}{\Gamma(a+b+2j+1)}P^{(a,b)}_j(1-2x) .
\label{}
\end{equation}
Let $ (p_j,p_j)_2 :=\int^1_0 (p_j(x))^2x^a(1-x)^bdx $, and define
\begin{equation}
   \tilde{K}^J_N(x,y) = \frac{(w_2(x)w_2(y))^{1/2}}{(p_{N-1},p_{N-1})_2}
                        \frac{p_{N}(x)p_{N-1}(y)-p_{N-1}(x)p_{N}(y)}{x-y} .
\label{Jkernel}
\end{equation}
It is well known, and easy to derive (see e.g. \cite{rmt_Fo}), that with 
$  \tilde{K}^J_{N,(t,1)} $ denoting the integral operator on $ (t,1) $ with kernel
(\ref{Jkernel}),
\begin{equation}
   E^J_N(t;a,b;\xi) = \det(1-\xi\tilde{K}^J_{N,(t,1)}) .
\label{JUE_Fdet}
\end{equation}
It was in this form that
\begin{equation} 
  t(t-1)\frac{d}{dt}\log E^J_N
\label{JUE_sigma}
\end{equation} 
was first related to
the solution of a nonlinear equation by Tracy and Widom \cite{TW_1994}. The 
equation found was of third order. Subsequently Haine and Semengue \cite{HS_1999}
studied (\ref{EJUE}) itself in the case $ \xi=1 $ and found a different third
order equation for (\ref{JUE_sigma}). Upon subtracting the two equations they
obtained a second order second degree nonlinear equation which they identified 
as an example of the $ \sigma$-form of Painlev\'e VI (\ref{PVI_sigmaF}). The
study of Tracy and Widom proceeded via functional properties of quantities 
associated with the Fredholm determinant (\ref{JUE_Fdet}), while Haine and 
Semengue used the theory of the KP hierarchy and Virasoro constraints satisfied
by certain matrix integrals as introduced by Adler and van Moerbeke \cite{AvM_2001}.
A third approach to the problem was initiated by Borodin and Deift \cite{BD_2002}.
They combined Riemann-Hilbert theory with the method of isomonodromic deformation
of certain linear differential equations to obtain a characterisation of 
(\ref{JUE_sigma}) which allows for immediate identification with the parameters
in (\ref{PVI_sigmaF}) (this is not the case with \cite{HS_1999}). Explicitly 
it was shown that
\begin{equation} 
  \sigma(t) = -t(t-1)\frac{d}{dt}\log E^J_N(1-t;a,b;\xi)
  +v_1v_2t+\frac{1}{2}(-v_1v_2+v_3v_4) ,
\label{JUE_BDsigma}
\end{equation}
satisfies (\ref{PVI_sigmaF}) with
\begin{equation}
   v_1=v_2=N+\frac{a+b}{2}, \quad v_3=\frac{a+b}{2}, \quad v_4=\frac{a-b}{2} .
\label{}
\end{equation}    
Furthermore, it is required that as $ t\to 0 $ 
\begin{align}
  \frac{d}{dt}\log E^J_N(1-t;a,b;\xi) 
 & \sim -\xi \tilde{K}^J_N(1-t,1-t)
 \\
 & = -\xi C_N(a,b) (1-t)^b ,
\label{}
\end{align} 
where 
\begin{equation}
  C_N(a,b) = \frac{\Gamma(a+b+N+1)\Gamma(b+N+1)}{\Gamma(N)\Gamma(a+N)\Gamma(b+1)\Gamma(b+2)} ,
\label{}
\end{equation}
thus providing the boundary condition to be satisfied by (\ref{JUE_BDsigma}).

According to (\ref{JUE_BDsigma}) 
\begin{equation} 
  E^J_N(1-t;a,b;\xi) = \exp\int^1_{1-t} \frac{ds}{s(1-s)}
  \left(\sigma(s)-v_1v_2s-\frac{1}{2}(-v_1v_2+v_3v_4)\right) .
\label{JUE_BDtau}
\end{equation}
From the characterisation of $ \sigma(t) $ as a $ \sigma$\PVI transcendent with
a specific boundary condition, the power series solution of (\ref{JUE_BDtau}) 
about $ t=1 $ can be readily computed \cite{rmt_Fo},
\begin{multline}
  E^J_N(1-t;a,b;\xi) = 1 \\
  -\xi\frac{C_N(a,b)}{b+1}(1-t)^{b+1}\left\{ 
            1-\frac{(b+1)(2N^2+2(a+b)N-2-2b+ab)}{(b+2)^2}(1-t)+{\rm O}((1-t)^2) \right\} \\
   +\xi^2\frac{C^2_N(a,b)(N-1)(N+b+1)(N+a-1)(N+a+b+1)}{(b+2)^2(b^2+4b+3)^2}(1-t)^{2b+4}
    \{1+{\rm O}(1-t)\} .
\label{JUE_BC}
\end{multline}
Furthermore, one can anticipate from (\ref{JUE_Fdet}) that the leading term in 
$ 1-t $ accompanying the power $ \xi^k $ will be proportional to 
$ (1-t)^{kb+k^2} $, as is consistent with (\ref{JUE_BC}).

It is well known (see e.g. \cite{rmt_Fo}) that after the change of variables 
$ x_j=\cos^2\theta_j/2, 0\leq \theta_j < \pi $ and with $ a,b = \pm 1/2 $
the eigenvalue p.d.f for the Jacobi ensemble as specified by (\ref{Epdf}) and
(\ref{classicalW}) becomes identical to the eigenvalue p.d.f. for matrices
from the classical groups $ O^{\pm}(N), Sp(2N) $ chosen with Haar (uniform) 
measure. As a consequence (in an obvious notation)
\begin{align}
   E^{O^-(2N+1)}((0,\phi);\xi)
   & = \left. E^J_N(\cos^2\phi/2;\xi)\right|_{{a=1/2}\atop{b=-1/2}}
   \nonumber \\
   E^{O^+(2N+1)}((0,\phi);\xi)
   & = \left. E^J_N(\cos^2\phi/2;\xi)\right|_{{a=-1/2}\atop{b=1/2}} .
\end{align}
Analogous to the expansion (\ref{JUE_BC}), the $ \sigma$\PVI evaluation 
(\ref{JUE_BDtau}) can then be used to deduce the expansions \cite{rmt_Fo}
\begin{multline}
  E^{O^-(2N+1)}((0,x);\xi) = 1-\tilde{c}x+\frac{4N^2-1}{36}\tilde{c}x^3 \\
  -\frac{48N^4-40N^2+7}{3600}\tilde{c}x^5 
  +\frac{4N^4-5N^2+1}{2025}\tilde{c}^2x^6 \\
  +\frac{192N^6-336N^4+196N^2-31}{211680}\tilde{c}x^7
  -\frac{48N^6-112N^4+77N^2-13}{198450}\tilde{c}^2x^8 + {\rm O}(x^9) , 
\end{multline}
\begin{multline}
  E^{O^+(2N+1)}((0,x);\xi) = 1-\frac{4N^2-1}{36}\tilde{c}x^3 \\
  +\frac{(4N^2-1)(12N^2-7)}{3600}\tilde{c}x^5 
  -\frac{(4N^2-1)(48N^4-72N^2+31)}{211680}\tilde{c}x^7
  + {\rm O}(x^9) , 
\end{multline}
where $ \tilde{c} = 2N\xi/\pi $.

The method of \cite{TW_1994} was adopted in \cite{WF_2000} to relate
\begin{multline}
  E^{Cy}_N(s;\eta;\xi) 
  := \frac{1}{C}\left(\int^{\infty}_{-\infty}-\xi\int^{\infty}_s\right)dx_1 \cdots
                \left(\int^{\infty}_{-\infty}-\xi\int^{\infty}_s\right)dx_N \\
     \times \prod^N_{l=1}\frac{1}{(1+x^2_l)^{\eta}} \prod_{1\leq j<k \leq N}(x_k-x_j)^2 ,
\label{ECyUE}
\end{multline}
to $ \sigma$\PVI. Explicitly it was shown that 
\begin{equation} 
  \sigma(s) = (1+s^2)\frac{d}{ds}\log E^{Cy}_N(s;a+N;\xi),
\label{CyUE_FWsigma}
\end{equation}
satisfies the equation
\begin{multline}
  (1+s^2)^2(\sigma'')^2+4(1+s^2)(\sigma')^3-8s\sigma(\sigma')^2+4\sigma^2(\sigma'-a^2) \\
  +8a^2s\sigma\sigma'+4[N(N+2a)-a^2s^2](\sigma')^2 = 0 .
\label{CyUE_sigmaF}
\end{multline}
As noted in \cite{FW_2004}, the relationship between (\ref{CyUE_sigmaF}) and 
(\ref{PVI_sigmaF}) can be seen by changing variables 
\begin{equation}
   t \mapsto \frac{is+1}{2}, \qquad \sigma_{\rm VI}(t) \mapsto \frac{i}{2}h(s) ,
\end{equation}
in the latter so that it reads
\begin{equation}
   h'\left((1+s^2)h''\right)^2+4\left(h'(h-sh')-iv_1v_2v_3v_4\right)^2
   +4\prod^4_{k=1}(h'+v^2_k) = 0 .
\label{PVI_sigmaFi}
\end{equation}
With 
\begin{equation}
   h=\sigma-a^2s, \quad v_1=-a, \quad v_2=0, \quad v_3=N+a, \quad v_4=a ,
\end{equation}
(\ref{PVI_sigmaFi}) reduces to (\ref{CyUE_sigmaF}).

The interest in (\ref{ECyUE}) in random matrix theory comes about by making a
stereographic projection from the real line to the unit circle by the change 
of variables  $ x = \tan\theta/2 $. With $ z = e^{i\theta} $ this shows 
\begin{multline}
   \prod^N_{l=1}\frac{1}{(1+x^2_l)^{N+\eta}}
   \prod_{1\leq j<k \leq N}(x_k-x_j)^2 dx_1 \cdots dx_N \\
   = 2^{-N(N+2\eta)}\prod^N_{l=1}|1+z_l|^{2\eta}
                    \prod_{1\leq j<k \leq N}|z_k-z_j|^2 d\theta_1 \cdots d\theta_N .
\label{Cy2cJ}
\end{multline}
In the case $ \eta = 0 $ the measure on the right hand side of (\ref{Cy2cJ})
corresponds to the eigenvalue p.d.f. for $ N\times N $ random unitary matrices
chosen with Haar (uniform) measure. For general $ \eta \in \NN $ it corresponds
to this same ensemble conditional so that there is an eigenvalue of degeneracy
$ \eta $ at $ \theta=\pi $.

With $ E^{U(N)}_N((\phi_1,\phi_2);\xi) $ denoting the generating function for 
the probability that the interval $ (\phi_1,\phi_2) $ contains exactly $ n $
eigenvalues, it follows that
\begin{equation}
  E^{U(N)}_N((0,2x);\xi) = \exp\left(-\int^x_0 h(\cot\phi)d\phi\right) ,
\end{equation} 
where $ h(s) $ satisfies (\ref{PVI_sigmaFi}) with $ v_1=v_2=v_3=0, v_4=N $.
Since for $ x \to 0 $, $ E^{U(N)}_N((0,2x);\xi) \sim 1-\xi xN/\pi $, we seek 
the solution of (\ref{PVI_sigmaFi}) subject to the boundary condition 
$ h(s) \sim c $, $ c:= \xi N/\pi $. This allows the power series expansion
\begin{multline}
   E^{U(N)}_N((0,2x);\xi) = 1-cx+\frac{N^2-1}{36}c^2x^4
  -\frac{(N^2-1)(2N^2-3)}{1350}c^2x^6 \\
  +\frac{(N^2-1)(N^2-2)(3N^2-5)}{52920}c^2x^8 
  -\frac{(N^2-4)(N^2-1)^2}{291600}c^3x^9 + {\rm O}(x^{10}) ,
\label{}
\end{multline}
to be computed \cite{rmt_Fo}. This expansion was first computed in \cite{TW_1994}
using the characterisation of $ E^{U(N)}_N $ in terms of a third order nonlinear
differential equation.

We have given reference to three distinct approaches which relate (\ref{EJUE})
to nonlinear differential equations with the Painlev\'e property. There is a 
fourth approach, which is due to the present authors \cite{FW_2001a}, \cite{FW_2002a}
\cite{FW_2004}, and involves applying Okamoto's theory of the Hamiltonian systems
approach to \PVI \cite{Ok_1987a}. This approach has the advantage of allowing 
generalisations of the generating functions (\ref{EJUE}), (\ref{ECyUE}) and 
$ E^{U(N)}_N((\phi_1,\phi_2);\xi) $ to be related to $ \sigma$\PVI. Consider 
for definiteness the latter. In \cite{FW_2004} the more general quantity
\begin{multline}
   A_N(t;\omega_1,\omega_2,\mu;\xi^*)
    := {1\over N!} 
       \left(\int^{\pi}_{-\pi}-\xi^*\int^{\pi}_{\pi-\phi}\right) {d\theta_1 \over 2\pi}
\ldots \left(\int^{\pi}_{-\pi}-\xi^*\int^{\pi}_{\pi-\phi}\right) {d\theta_N \over 2\pi} \\
      \times\prod^{N}_{l=1} z^{-i\omega_2}_l |1+z_l|^{2\omega_1} |1+tz_l|^{2\mu}
      \prod_{1 \leq j < k \leq N} |z_j-z_k|^2 ,
\label{A}
\end{multline}
where $ t = e^{i\phi}, \phi \in [0,2\pi) $, $ \xi^* \in \CC $ and the 
parameters $ \omega_1,\omega_2,\mu \in \CC $, $ \omega = \omega_1+i\omega_2 $,
are restricted with $ \Re(2\omega_1), \Re(2\mu) > -1 $, $ N \in \ZZ_{\geq 0} $. 
The independent variable $ t $, whilst originally defined on the unit circle 
$ |t| = 1 $ with a real angle $ \phi $, can be considered as a complex variable
which is analytically continued into the cut complex $t$-plane.
The case $\omega_2 = \mu = 0$, $\omega_1=\eta$ of (\ref{A}) gives the
generating function for the probability of $k$ eigenvalues in $(0,\phi)$
for the ensemble specified by the right hand side of (\ref{Cy2cJ}).
Define
\begin{align}
   M_N(a,b)
  & := \int^{1/2}_{-1/2}dx_1 \ldots \int^{1/2}_{-1/2}dx_N
       \prod^{N}_{l=1}z^{(a-b)/2}_l |1+z_l|^{a+b}
       \prod_{1 \leq j < k \leq N} |z_j-z_k|^2 \\
  &  = \prod^{N-1}_{j=0}\frac{\Gamma(a+b+j+1)\Gamma(j+2)}{\Gamma(a+j+1)\Gamma(b+j+1)},
     \qquad z_l = e^{2\pi i x_l}, l=1,\ldots,N ,
\label{morrisI}
\end{align}
set $ \xi^*=1-(1-\xi)e^{-\pi i\mu} $ and denote by $ \tilde{h}(s) $ a solution
of (\ref{PVI_sigmaFi}) with
\begin{equation}
   v_1=-\mu-\omega_1, \quad v_2=i\omega_2, \quad 
   v_3=N+\mu+\omega_1, \quad v_4=-\mu+\omega_1 .
\label{AN_paramCyUE}
\end{equation}
It was shown in \cite{FW_2004} that 
\begin{multline}
   A_N(t;\omega_1,\omega_2,\mu;\xi^*)
    = \frac{M_N(\mu+\bar{\omega},\mu+\omega)}{M_N(0,0)} \\ \times
      \exp\left( -\frac{1}{2}\int^{\phi}_0\left[ \tilde{h}(\cot\frac{\theta}{2})
        +\omega_2(N+\omega_1-\mu)+(\mu+\omega_1)^2\cot\frac{\theta}{2} \right]d\theta
          \right) .
\label{}
\end{multline}
We remark that according to the theory in the sentence including (\ref{PVI_sigmaFi})
\begin{equation}
   \tilde{h}(\cot\frac{\theta}{2}) = -2i\sigma_{\rm VI}(\frac{1}{1+e^{i\theta}}) .
\label{xfm}
\end{equation}

A generalisation of the Jacobi ensemble (\ref{EJUE}) may also be
made in a similar fashion, by the introduction of the additional
factor $ \prod_{l=1}^N|s - x_l|^\mu $ in the integrand. In the case
$\xi=0$ it is shown in \cite{FW_2004} that this generalised ensemble can be
related to (\ref{A}), and consequently that
\begin{multline}
   \sigma_{\rm VI}(t) 
  = \left( e_2[-v_1,-v_2,v_3]+N\mu \right)t-\frac{1}{2}e_2[-v_1,-v_2,v_3,v_4]-N\mu 
  \\
    + t(t-1){d \over dt}\log A_N(t;\omega_1,\omega_2,\mu;\xi^*=0),
\label{SSE_PVIsoln}
\end{multline}
is a solution of (\ref{PVI_sigmaF}) with the parameters
\begin{equation}
   v_1 = {\displaystyle N+\omega-\mu \over \displaystyle 2}, \quad
   v_2 = \bar{\omega}+{\displaystyle N+\omega+\mu \over \displaystyle 2}, \quad
   v_3 = {\displaystyle N-\omega+\mu \over \displaystyle 2}, \quad
   v_4 = -\mu-{\displaystyle N+\omega+\mu \over \displaystyle 2} .
\label{AN_paramJUE}
\end{equation}
Moreover a derivation of the transformation implied by (\ref{xfm}) and
(\ref{SSE_PVIsoln}) was given which tells us that the condition $\xi^*=0$ in
(\ref{SSE_PVIsoln}) can be relaxed. 

The above mentioned generalisation of the Jacobi ensemble (\ref{EJUE}) has 
also revealed different Painlev\'e connections to those found in
\cite{FW_2004}. In the mid 1990's, with $\xi=0$, this was studied from the
viewpoint of orthogonal polynomial theory by Magnus \cite{Ma_1995a}.
In that study an auxiliary quantity occurring in the theory was shown
to satisfy the \PVI equation (\ref{PVI_y}) for appropriate parameters. We 
remark that further development of orthogonal polynomial theory in relation 
to (\ref{A}) \cite{FW_2004a,FW_2004b} has been shown to relate to 
B\"acklund transformations in the Hamiltonian theory of \PVI, and to the 
so called discrete \dPV equations.

In the case $ \omega_2=\mu=0, \omega_1=\eta $ (\ref{A}) is the generating function
for the probability that there are exactly $ n $ eigenvalues in 
$ (\pi-\phi,\pi) $ for the eigenvalue p.d.f. on the right hand side of (\ref{A}).
In this case expanding (\ref{JUE_Fdet}) allow us to specify the $ s \to \infty $
boundary condition which must be satisfied by $ \tilde{h}(s) $. However, for generic
parameters the boundary conditions were not given in \cite{FW_2004}. Consideration 
of this latter problem is our concern for the remainder of the paper.
A similar situation was recently rectified in \cite{FW_2006} in relation to
a class of multidimensional integral solutions of $ \sigma$\PV \cite{FW_2002a}.
The approach taken was to write the multidimensional integral as a determinant
(involving the confluent hypergeometric function), and to expand the entries of
the determinant. An analogous strategy suffices in relation to (\ref{A}),
now with the entries of the determinant given in terms of the Gauss hypergeometric
function. This is done in Section \ref{section_BCSSE}.

The form of the general expansion of a solution of the $ \sigma$\PVI equation
(\ref{PVI_sigmaF}), or equivalently its $ \tau$ function about its fixed 
singularities $ t=0,1,\infty $ has been given by Jimbo \cite{Ji_1982}. The 
exponents occurring therein are given in terms of the monodromy data associated 
with the isomonodromy deformation formulation of \PVI. Aspects of this theory
and the results of Jimbo are revised in Section \ref{sectionPVI_isoMD}. 

One of the features of the general expansions is that they consist of two branches,
whereas the expansion of (\ref{A}) only exhibits a single branch. We show in 
Section \ref{sectionMData_SSE} that this is entirely consistent with the monodromy 
data associated 
with $ \sigma_{\rm VI} $; it is such that the coefficient in front of one of the
branches vanishes identically.

\section{Boundary Conditions}\label{section_BCSSE}
\setcounter{equation}{0}

Consider the multidimensional integral (\ref{A}). According to 
(\ref{xfm}) (with the condition $ \xi^*=0 $ relaxed), we know how to
relate its logarithmic derivative to a solution of the $\sigma$\PVI
equation (\ref{PVI_sigmaF}). However this property cannot be used to characterise
the multidimensional integral unless an appropriate boundary condition 
is specified. As remarked in the second paragraph below (\ref{AN_paramJUE})
above, in \cite{FW_2004} the interpretation of (\ref{A}) for some special
parameters as the generating function of a gap probability provided
the boundary condition in those cases. However for general parameters
no boundary condition was presented.

Taking a different approach, namely the expansion of the elements in the 
determinant form of (\ref{A}), the sought boundary conditions can be deduced as 
presented in the following result.

\begin{proposition}
For generic values of $ \mu, \omega, \bar{\omega} $ 
the spectral average $ A_N $ has the following expansions. About $ t=0 $ 
subject to $ \mu-\bar{\omega} \notin \ZZ $ we have
\begin{multline}
  t^{N\mu}A_N \mathop{\sim}\limits_{t \to 0}
  \left\{ 1+\xi^*{e^{-\pi i(\mu-\bar{\omega})} \over 2i\sin\pi(\mu-\bar{\omega})} \right\}^N
  \prod^{N-1}_{k=0}{k!\Gamma(2\omega_1+k+1) \over 
                      \Gamma(1+k+\mu+\omega)\Gamma(1+k-\mu+\bar{\omega})} \\
  \times\Bigg\{ 1 + {2N\mu(\mu+\omega) \over N-\mu+\bar{\omega}}t \\
          - {\xi^*e^{-\pi i(\mu-\bar{\omega})} \over 
             2i\sin\pi(\mu-\bar{\omega})+\xi^*e^{-\pi i(\mu-\bar{\omega})}} \\
  \times
        {\Gamma(1+\mu+\omega)\Gamma(1+\mu-\bar{\omega})\Gamma(1+2\mu)\Gamma(N-\mu+\bar{\omega}) 
         \over 
         \Gamma(N)\Gamma(N+\mu+\bar{\omega})\Gamma(N+2\omega_1)\Gamma^2(2-N+\mu-\bar{\omega})}
         t^{1-N+\mu-\bar{\omega}} \Bigg\} .
\label{AN_exp_0}
\end{multline}
About $ t=1 $ subject to $ 2\mu+2\omega_1 \notin \ZZ $ we have
\begin{multline}
  A_N \mathop{\sim}\limits_{t \to 1}
  \prod^{N-1}_{k=0}{k!\Gamma(2\mu+2\omega_1+k+1) \over 
                      \Gamma(1+k+\mu+\omega)\Gamma(1+k+\mu+\bar{\omega})}
  \Bigg\{ 1 + {N\mu(\bar{\omega}-\omega) \over 2\mu+2\omega_1}(1-t) \\
  + \frac{(-1)^{N+1}}{\sin\pi(2\mu+2\omega_1)}
    \left(\xi^*{e^{-\pi i(\mu-\bar{\omega})} \over 2i}
            +{\sin\pi 2\mu\sin\pi(\mu+\omega) \over \sin\pi(2\mu+2\omega_1)} \right) \\
  \times
  {\Gamma(1+2\mu)\Gamma(1+2\omega_1)\Gamma(1+\mu+\omega)\Gamma(1+\mu+\bar{\omega})
         \over 
   \Gamma^2(2\mu+2\omega_1+2)\Gamma(2\mu+2\omega_1+1)\Gamma(N)\Gamma(-N-2\mu-2\omega_1) }
         (1-t)^{1+2\mu+2\omega_1} \Bigg\} .
\label{AN_exp_1}
\end{multline}
And about $ t=\infty $ subject to $ \mu-\omega \notin \ZZ $ we have
\begin{multline}
  t^{-\mu N}A_N \mathop{\sim}\limits_{t \to \infty}
  \left\{ -\frac{e^{-\pi i2\mu}}{\sin\pi(\mu-\omega)}
           \left(\sin\pi(\mu+\omega)+\xi^*\frac{e^{-\pi i(\mu+\omega)}}{2i}
           \right) \right\}^N \\
  \times \prod^{N-1}_{k=0}\frac{k!\Gamma(2\omega_1+k+1)}
                     {\Gamma(1+k+\mu+\bar{\omega})\Gamma(1+k-\mu+\omega)} \\
  \times\Bigg\{ 1 + \frac{2\mu N(\mu+\bar{\omega})}{N-\mu+\omega}\frac{1}{t} \\
          + e^{\pi i(\mu-\omega)}\left( \frac{2i\sin\pi2\mu+\xi^*e^{-\pi i2\mu}}
              {2i\sin\pi(\mu+\omega)+\xi^*e^{-\pi i(\mu+\omega)}} \right) \\
  \times{\Gamma(1+2\mu)\Gamma(1+\mu+\bar{\omega})\Gamma(1+\mu-\omega)\Gamma(N-\mu+\omega) 
         \over 
         \Gamma(N)\Gamma(N+\mu+\omega)\Gamma(N+2\omega_1)\Gamma^2(2-N+\mu-\omega)}
         t^{-1+N-\mu+\omega} \Bigg\}
\label{AN_exp_infty}
\end{multline}
Any one of these boundary conditions suffices to uniquely define a solution to the 
ordinary differential equation (\ref{PVI_sigmaF}) under the above generic restrictions
on the parameters. 
\end{proposition}
\begin{proof}
Using Heine's identity we can deduce from (\ref{A}) that
\begin{multline}
   A_N(t;\omega_1,\omega_2,\mu;\xi^*) \\
     = t^{-N\mu}\det\left[
       \left(\int^{\pi}_{-\pi}-\xi^*\int^{\pi}_{\pi-\phi}\right) {d\theta \over 2\pi}
        z^{-\mu-\omega+j-k} (1+z)^{2\omega_1} (1+tz)^{2\mu} \right]_{0 \leq j,k \leq N-1} .
\label{SSE}
\end{multline}
This is a determinant of a Toeplitz matrix whose
symbol, or more specifically weight $ w(z) $, is one defined on the unit circle 
$ |z|=1 $ by
\begin{equation}
   w(z) = t^{-\mu}z^{-\mu-\omega} (1+z)^{2\omega_1} (1+tz)^{2\mu}
   \begin{cases}
    1, \quad \theta \in (-\pi,\pi-\phi) \cr
    1-\xi^*, \quad \theta \in (\pi-\phi,\pi) \cr
   \end{cases} ,
\end{equation}
with Fourier components $ w_k $ such that $ w(z) = \sum^{\infty}_{k=-\infty}w_k z^k $.
Considering the expansion about $ t=0 $ first we note that the Toeplitz matrix element 
can be evaluated as
\begin{multline}
   t^{\mu}w_{n} =
   \left\{ 1+\xi^*{e^{-\pi i(n+\mu-\bar{\omega})} \over
           2i\sin\pi(n+\mu-\bar{\omega})}
   \right\}
   {\Gamma(2\omega_1+1) \over
                    \Gamma(1+n+\mu+\omega)\Gamma(1-n-\mu+\bar{\omega})} \\
   \times{}_2F_1(-2\mu,-n-\mu-\omega;1-n-\mu+\bar{\omega};t) \\
   -\xi^*{e^{-\pi i(n+\mu-\bar{\omega})} \over 2i\sin\pi(n+\mu-\bar{\omega})}
       {\Gamma(2\mu+1) \over
                    \Gamma(1+n+\mu-\bar{\omega})\Gamma(1-n+\mu+\bar{\omega})} \\
           \times t^{n+\mu-\bar{\omega}}
   {}_2F_1(-2\omega_1,n-\mu-\bar{\omega};1+n+\mu-\bar{\omega};t) ,
\label{VI_toepM1}
\end{multline}
under the condition that $ n+\mu-\bar{\omega} \notin \ZZ $,
which is a form suitable for the development of an expansion about $ t=0 $.
In the course of the derivation we had to invoke tighter constraints on the 
parameters, namely $ |t|=1, t \neq 1 $ and $ \Re(2\mu) \geq 0, \Re(2\omega_1) \geq 0 $,
but these can be relaxed by analytic continuation arguments.
The structure is $ t^{\mu}w_n = a_n(t) + t^{n+\mu-\bar{\omega}}b_n(t) $ where
$ a_n(t), b_n(t) $ are analytic about $ t=0 $. 
One can expand the Toeplitz determinant about $ t=0 $ retaining only the leading
order terms from both the analytic and non-analytic contributions and the
resulting formula is 
\begin{multline}
   t^{\mu N}A_N(t) \mathop{\sim}\limits_{t \to 0}
   \det(a_{j-k}(0)+ta'_{j-k}(0))_{j,k=0,\ldots,N-1} \\
   + (-1)^{N-1}b_{-(N-1)}\det(a_{j-k+1}(0))_{j,k=0,\ldots,N-2}t^{\mu-\bar{\omega}-(N-1)} .
\label{Det_T=0Exp}
\end{multline}
Using the determinant identity \cite{Nd_2004} 
\begin{equation}
   \det\left(\frac{\Gamma(d+k-j)}{\Gamma(c+k-j)}\right)_{0 \leq j,k \leq n-1}
   = \prod^{n-1}_{j=0} j!\frac{\Gamma(1+d-c)}{\Gamma(1+d-c-j)}
                         \frac{\Gamma(d-n+1+j)}{\Gamma(c+j)} ,
\end{equation}
or its slightly more general form
\begin{multline}
   \det\left(\frac{\Gamma(z_k+b-j)}{\Gamma(z_k-j)}\right)_{0 \leq j,k \leq n-1} \\
   = \prod_{0 \leq j<k \leq n-1}(z_k-z_j)
     \prod^{n-1}_{j=0} (-1)^j(-b)_{j}\frac{\Gamma(z_j+b-n+1)}{\Gamma(z_j)} ,
\end{multline}
for an arbitrary sequence $ \{z_j\}^{n-1}_{j=0} $ one can evaluate the determinants 
appearing in (\ref{Det_T=0Exp}). The result is (\ref{AN_exp_0}).

To develop the expansion about $ t=1 $ the appropriate expression for the 
Toeplitz matrix element is
\begin{multline}
   t^{\bar{\omega}-n}w_{n} =
   {\Gamma(2\mu+2\omega_1+1) \over
                    \Gamma(1+n+\mu+\omega)\Gamma(1-n+\mu+\bar{\omega})} \\
   \times{}_2F_1(-2\omega_1,n-\mu-\bar{\omega};-2\mu-2\omega_1;1-t) \\
   +\frac{1}{\pi}
\left(\xi^*{e^{-\pi i(n+\mu-\bar{\omega})} \over 2i}
            +{\sin\pi 2\mu\sin\pi(n+\mu+\omega) \over \sin\pi(2\mu+2\omega_1)} \right)
       {\Gamma(1+2\mu)\Gamma(1+2\omega_1) \over \Gamma(2+2\mu+2\omega_1)} \\
           \times (1-t)^{1+2\mu+2\omega_1}
   {}_2F_1(1+2\mu,n+1+\mu+\omega;2+2\mu+2\omega_1;1-t) ,
\label{VI_toepM2}
\end{multline}
subject to $ 1+2\mu+2\omega_1 \notin \ZZ $.
The structure is now 
$ t^{\bar{\omega}-n}w_n = a_n(1-t) + (1-t)^{1+2\mu+2\omega_1}b_n(1-t) $ 
where $ a_n(1-t), b_n(1-t) $ are analytic about $ t=1 $. Again one
can expand the Toeplitz determinant about $ t=1 $ retaining only the leading
order terms from both the analytic and non-analytic contributions and use the above
determinant identities to arrive at (\ref{AN_exp_1}).

Lastly the expansion about $ t=\infty $ is computed using the Toeplitz matrix 
element in the form
\begin{multline}
   t^{\mu}w_{n} =
   \frac{e^{-\pi i2\mu}}{\sin\pi(n-\mu+\omega)}
    \left( \sin\pi(n+\mu+\omega)+\xi^*\frac{e^{-\pi i(n+\mu+\omega)}}{2i} \right) \\
   \times \frac{\Gamma(2\omega_1+1)}
                   {\Gamma(1+n-\mu+\omega)\Gamma(1-n+\mu+\bar{\omega})} \\
   \times t^{2\mu}{}_2F_1(-2\mu,n-\mu-\bar{\omega};1+n-\mu+\omega;1/t) \\
   -\frac{e^{-\pi i(n+\mu+\omega)}}{\sin\pi(n-\mu+\omega)}
    \left( \sin\pi2\mu+\xi^*\frac{e^{-\pi i2\mu}}{2i} \right)
    \frac{\Gamma(2\mu+1)}{\Gamma(1+n+\mu+\omega)\Gamma(1-n+\mu-\omega)} \\
           \times t^{n+\mu+\omega}
   {}_2F_1(-2\omega_1,-n-\mu-\omega;1-n+\mu-\omega;1/t) ,
\label{VI_toepM3}
\end{multline}
valid for $ n-\mu+\omega \notin \ZZ $.
The structure is $ t^{-\mu}w_n = a_n(1/t) + t^{n-\mu+\omega}b_n(1/t) $ where
$ a_n(1/t), b_n(1/t) $ are analytic about $ t=\infty $. Repeating the procedure adopted 
about the other singularities one arrives at (\ref{AN_exp_infty}).
\end{proof}

\begin{remark}
We do not attempt to treat the degenerate cases where either
$ \mu-\bar{\omega} \in \ZZ $, $ 2\mu+2\omega_1 \in \ZZ $ or $ \mu-\omega \in \ZZ $
here as this results in confluent logarithms and other technical difficulties.
\end{remark}

\section{Isomonodromy Deformation Formulation for \PVI}\label{sectionPVI_isoMD}
\setcounter{equation}{0}

In this section we describe the isomonodromic deformation system which characterises
the general solution to the sixth Painlev\'e equation and outline a solution to 
the direct monodromy problem, that is explicit formulae for the monodromy data
associated with a particular solution of the $\sigma$-form for the sixth Painlev\'e 
equation.
The isomonodromic deformation system associated with the sixth Painlev\'e equation
is not uniquely determined and we shall see this arbitrariness arise in 
Section \ref{sectionMData_SSE} 
when applying the general theory to our random matrix theory.
  
Following the conventions and notations of \cite{Ji_1982}, \cite{Ki_1994} we consider 
the $2\times 2$ linear matrix ODE for $ \Psi(\lambda;t) $ with four regular 
singularities in the $\lambda$-plane chosen to be $ \nu = 0,t,1,\infty $
\begin{align}
   {\partial \over \partial\lambda}\Psi 
  & = \left( {A_{0} \over \lambda}+{A_{1} \over \lambda-1}+{A_{t} \over \lambda-t}
      \right) \Psi ,
  \label{PVI_linear:a} \\
    {\partial \over \partial t}\Psi 
  & = -{A_{t} \over \lambda-t} \Psi .
  \label{PVI_linear:b}
\end{align}
It is taken that the residue matrices $ A_{\nu}(t) $ satisfy
\begin{equation}
   A_{0}+A_{t}+A_{1}= -A_{\infty} = -{\theta_{\infty} \over 2}\sigma_3, \quad
   \sigma_3 := \left( \begin{array}{cc} 1 & 0 \\ 0 & -1 \end{array} \right), \quad
   \theta_{\infty} \in \CC\backslash\ZZ,
\label{PVI_Ainfty}
\end{equation}
and given the freedom to choose either $ {\rm tr} A_{\nu} $ or one of the eigenvalues
we follow the convention of Jimbo \cite{Ji_1982}
\begin{equation}
   {\rm tr} A_{\nu} = 0, \quad \det A_{\nu} = -\frac{1}{4} \theta^2_{\nu},
    \; \nu = 0,t,1,\infty ,
\label{PVI_Aint}
\end{equation}
defining the formal exponents of monodromy $ \theta_{\nu} $.
The residue matrices satisfy the {\it Schlesinger system} of equations 
\begin{align}
   {d \over dt}A_{0}
  & = -\frac{[A_{0},A_{t}]}{t} ,
  \label{PVI_Schles:a} \\
   {d \over dt}A_{t}
  & =  \frac{[A_{0},A_{t}]}{t}+\frac{[A_{t},A_{1}]}{1-t} , 
  \label{PVI_Schles:b} \\
   {d \over dt}A_{1}
  & = -\frac{[A_{t},A_{1}]}{1-t} , 
  \label{PVI_Schles:c} \\
   {d \over dt}A_{\infty}
  & = 0 , 
  \label{PVI_Schles:d}
\end{align}
as a consequence of the compatibility of (\ref{PVI_linear:a}) and (\ref{PVI_linear:b}).

The {\it $\tau$-function} for \PVI is defined by
\begin{equation}
   \frac{d}{dt}\log\tau
   = {\rm Tr}\left({A_{0} \over t}+{A_{1} \over t-1} \right)A_{t} 
\label{PVI_tau}
\end{equation}
and the {\it $\sigma$-function} 
\begin{equation}
   \zeta(t) =
   t(t-1)\frac{d}{dt}\log\tau
   +\frac{1}{4}(\theta^2_{t}-\theta^2_{\infty})t
   -\frac{1}{8}(\theta^2_{t}+\theta^2_{0}-\theta^2_{\infty}-\theta^2_{1})
\label{PVI_sigma}
\end{equation}
which satisfies the second-order second degree differential equation
\begin{multline}
 \frac{d}{dt}\zeta\left(t(t-1)\frac{d^2}{dt^2}\zeta\right)^2 \\
 +\left[ 2\frac{d}{dt}\zeta\left(t\frac{d}{dt}\zeta-\zeta\right)
         -\left(\frac{d}{dt}\zeta\right)^2
         -\frac{1}{16}(\theta^2_{t}-\theta^2_{\infty})(\theta^2_{0}-\theta^2_{1})  
  \right]^2 \\
 =\left(\frac{d}{dt}\zeta+\frac{1}{4}(\theta_{t}+\theta_{\infty})^2\right)
  \left(\frac{d}{dt}\zeta+\frac{1}{4}(\theta_{t}-\theta_{\infty})^2\right) \\
  \times
  \left(\frac{d}{dt}\zeta+\frac{1}{4}(\theta_{0}+\theta_{1})^2\right)
  \left(\frac{d}{dt}\zeta+\frac{1}{4}(\theta_{0}-\theta_{1})^2\right) ,
\label{PVI_SF}
\end{multline}
(cf.~(\ref{PVI_sigmaF})).

Furthermore we suppose ({\it Assumption 1}) that the matrices $ A_{\nu} $ are 
diagonalisable, i.e. that there exists nonsingular $ R_{\nu} \in SL(2,\CC) $ such that
\begin{equation}
   R^{-1}_{\nu}A_{\nu}R_{\nu} = \frac{1}{2}\theta_{\nu}\sigma_3,
   \; \theta_{\nu} \in \CC\backslash\ZZ.
\label{PVI_Adiag}
\end{equation}
In the neighbourhood of a regular singularity $ \Psi(\lambda) $ can be expanded 
locally as
\begin{equation}
   \Psi = \sum^{\infty}_{m=0} \Psi_{m\nu}
            (\lambda-\nu)^{m+{\theta_{\nu} \over 2}\sigma_3} C_{\nu},
\label{PVI_regE}
\end{equation}
for $ \nu = 0,t,1 $ and for $ \lambda = \infty $ in the form
\begin{equation}
   \Psi = \left(I+ \sum^{\infty}_{m=1} \Psi_{m\infty}\lambda^{-m} \right)
            \lambda^{-{\theta_{\infty} \over 2}\sigma_3}.
\label{PVI_regEinfty}
\end{equation}
For such local series to exist we have to assume ({\it Assumption 2}) that the
eigenvalues of $ A_{\nu} $ are distinct modulo the non-zero integers, i.e.
that $ \pm\theta_{\nu} \notin \NN $. 

The matrices $ C_{\nu}, \nu=0,t,1 $ are the connection matrices and we are
taking the local solutions (\ref{PVI_regEinfty}) as our fundamental 
system of solutions, i.e. $ C_{\infty}=I $.
The monodromy matrices $ M_{\nu} (\nu = 0,t,1,\infty) $ are defined as
\begin{equation}
   \left.\Psi\right|_{\lambda=\nu+\delta e^{2\pi i}} 
 = \left.\Psi\right|_{\lambda=\nu+\delta}M_{\nu}, \qquad \nu=0,t,1,\infty ,
\label{PVI_Mdefn}
\end{equation}
and are given in terms of the monodromy exponents and connection matrices by
\begin{equation}
   M_{\nu} = C^{-1}_{\nu}e^{\pi i\theta_{\nu}\sigma_3}C_{\nu} .
\label{PVI_mon}
\end{equation}
The monodromy matrices together satisfy the cyclic relation
\begin{equation}
   M_{\infty}M_{1}M_{t}M_{0} = I ,
\end{equation}
according to the convention taken for the basis of loops displayed in Figure \ref{PVI_cyclic.fig},
which generate the fundamental group $ \pi(\CP\backslash\{0,t,1,\infty\},\lambda_0) $.
There is arbitrariness in the monodromy data in the sense that the replacement
$ C_{\nu} \mapsto D^{-1}_{\nu}C_{\nu} $ doesn't change the monodromy matrices 
provided that $ D_{\nu} $ commutes with the right-hand side of (\ref{PVI_Adiag}).
This implies that $ D_{\nu} $ is diagonal if $ \theta_{\nu} \neq 0 $.
This arbitrariness will manifest itself in the appearance of an arbitrary complex
number in the explicit parameterisation discussed latter.

\begin{figure}[H]
\begin{xy}{
   (0,0)*{};
   (60,20)="BASE";
   "BASE";"BASE"**\crv{(30,-40)&(50,-40)};?(0.25)*\dir{>};        
   "BASE";"BASE"**\crv{(50,-40)&(70,-40)};?(0.25)*\dir{>};        
   "BASE";"BASE"**\crv{(70,-40)&(90,-40)};?(0.25)*\dir{>};        
   "BASE";"BASE"**\crv{(90,-40)&(110,-40)};?(0.25)*\dir{>};        
   (47,-20)*+{\cdot\,0};(60,-20)*+{\cdot\,t};(74,-20)*+{\cdot\,1};(88,-20)*+{\cdot\,\infty};
   "BASE"+(0,3)*{\lambda_0}
}
\end{xy}
\vskip0.0cm
\caption{Monodromy representation of the fundamental group for 
$ \CP\backslash\{0,t,1,\infty\} $}\label{PVI_cyclic.fig}
\end{figure}

The isomonodromic principle states that the monodromy data 
$ MD:=\{ \theta_{\nu}, C_{\nu}, M_{\nu} | \nu = 0,t,1,\infty  \} $
are preserved under the deformations of $ t $.
The invariants of the monodromy data are defined
\begin{align}
   p_{\mu} &= 2\cos\pi\theta_{\mu} := {\rm Tr}M_{\mu} \quad
   \mu \in \{0,t,1,\infty\}
   \\
   p_{\mu\nu} &= 2\cos\pi\sigma_{\mu\nu} := {\rm Tr}M_{\mu}M_{\nu} \quad
   \mu,\nu \in \{0,t,1\} ,
\end{align}
in the sense that these do not contain any arbitrary constants.

In \cite{Ji_1982} Jimbo states the following conditions under which his results apply
\begin{align}
  & \theta_0, \theta_t, \theta_1, \theta_{\infty} \notin \mathbb{Z} ,
  \label{PVI_conditions:a}\\
  & 0 < {\Re}(\sigma_{0t}) < 1 ,
  \label{PVI_conditions:b}\\
  & \theta_0\pm\theta_t\pm\sigma_{0t},\quad 
    \theta_{\infty}\pm\theta_1\pm\sigma_{0t} \notin 2\mathbb{Z} .
  \label{PVI_conditions:c}
\end{align}

When $ \sigma_{0t} \neq 0 $ a parameterisation of the monodromy matrices was deduced
by Jimbo and is given in Lemma \ref{PVI_Mlemma}.
\vfill\eject

\begin{sideways}
\begin{minipage}[c][12cm][c]{20cm}{
\begin{lemma}[Jimbo\cite{Ji_1982}]\label{PVI_Mlemma}
Subject to the conditions (\ref{PVI_conditions:a}), (\ref{PVI_conditions:b}) and 
(\ref{PVI_conditions:c})
the monodromy matrices can be parameterised in the following way
\begin{equation}
  M_{\infty} = 
  \begin{pmatrix}
    e^{\pi i\theta_{\infty}} & 0 \cr
    0 & e^{-\pi i\theta_{\infty}}
  \end{pmatrix} ,
\label{PVI_Minfty}
\end{equation}
\begin{equation}
  M_{1} = \frac{1}{i\sin\pi\theta_{\infty}}
  \begin{pmatrix}
    \cos\pi\sigma-e^{-\pi i\theta_{\infty}}\cos\pi\theta_{1}
    &
    -2re^{-\pi i\theta_{\infty}}\sin\frac{\pi}{2}(\theta_{\infty}+\theta_{1}+\sigma)
                                  \sin\frac{\pi}{2}(\theta_{\infty}+\theta_{1}-\sigma)
  \cr
 2r^{-1}e^{\pi i\theta_{\infty}}\sin\frac{\pi}{2}(\theta_{\infty}-\theta_{1}+\sigma)
                                  \sin\frac{\pi}{2}(\theta_{\infty}-\theta_{1}-\sigma)
    &
    -\cos\pi\sigma+e^{\pi i\theta_{\infty}}\cos\pi\theta_{1}
  \end{pmatrix} ,
\label{PVI_M1}
\end{equation}
\begin{equation}
  CM_{t}C^{-1} = \frac{1}{i\sin\pi\sigma}
  \begin{pmatrix}
    e^{\pi i\sigma}\cos\pi\theta_{t}-\cos\pi\theta_{0}
    &
    -2se^{\pi i\sigma}\sin\frac{\pi}{2}(\theta_{0}+\theta_{t}-\sigma)
                                  \sin\frac{\pi}{2}(\theta_{0}-\theta_{t}+\sigma)
  \cr
 2s^{-1}e^{-\pi i\sigma}\sin\frac{\pi}{2}(\theta_{0}+\theta_{t}+\sigma)
                                  \sin\frac{\pi}{2}(\theta_{0}-\theta_{t}-\sigma)
    &
    -e^{-\pi i\sigma}\cos\pi\theta_{t}+\cos\pi\theta_{0}
  \end{pmatrix} ,
\label{PVI_Mt}
\end{equation}
\begin{equation}
  CM_{0}C^{-1} = \frac{1}{i\sin\pi\sigma}
  \begin{pmatrix}
    e^{\pi i\sigma}\cos\pi\theta_{0}-\cos\pi\theta_{t}
    &
     2s\sin\frac{\pi}{2}(\theta_{0}+\theta_{t}-\sigma)
                                  \sin\frac{\pi}{2}(\theta_{0}-\theta_{t}+\sigma)
  \cr
-2s^{-1}\sin\frac{\pi}{2}(\theta_{0}-\theta_{t}-\sigma)
                                  \sin\frac{\pi}{2}(\theta_{0}+\theta_{t}+\sigma)
    &
    -e^{-\pi i\sigma}\cos\pi\theta_{0}+\cos\pi\theta_{t}
  \end{pmatrix} ,
\label{PVI_M0}
\end{equation}
where
\begin{equation}
  C = 
  \begin{pmatrix}
    \sin\frac{\pi}{2}(\theta_{\infty}-\theta_{1}-\sigma)
    &
    r\sin\frac{\pi}{2}(\theta_{\infty}+\theta_{1}+\sigma)
    \cr
    r^{-1}\sin\frac{\pi}{2}(\theta_{\infty}-\theta_{1}+\sigma)
    &
    \sin\frac{\pi}{2}(\theta_{\infty}+\theta_{1}-\sigma)
  \end{pmatrix} .
\label{PVI_C}
\end{equation}
Here $ r $ is an arbitrary non-zero complex number, and the short hand notation
$ s:=s_{0t}, \sigma=\sigma_{0t} $ is used.
\end{lemma}}
\end{minipage} 
\end{sideways}
\vfill\eject

\begin{proof}
We begin by noting that (\ref{PVI_Minfty}) follows from our choice for the
fundamental system of solutions. To establish the other formula we require two 
preliminary results. Let us make the abbreviations 
$ c_{\nu}:=\cos\pi\theta_{\nu}, s_{\nu}:=\sin\pi\theta_{\nu},
  c_{\sigma}:=\cos\pi\sigma_{0t}, s_{\sigma}:=\sin\pi\sigma_{0t},
  \epsilon_{\nu}:=c_{\nu}+is_{\nu} $ and
\begin{equation}
  \mfS(\vartheta):=\sin\frac{\pi}{2}\vartheta ,
\end{equation}
and note some trigonometric identities which will be useful in this and the
ensuing proofs
\begin{align}
 2\mfS(\theta_{\infty}-\theta_1+\sigma)\mfS(\theta_{\infty}+\theta_1+\sigma)
 & = c_1-c_{\infty}c_{\sigma}+s_{\infty}s_{\sigma} ,
 \label{trigId:a} \\
 2\mfS(\theta_{\infty}+\theta_1-\sigma)\mfS(\theta_{\infty}-\theta_1-\sigma)
 & = c_1-c_{\infty}c_{\sigma}-s_{\infty}s_{\sigma} ,
 \label{trigId:b} \\
 2\mfS(\theta_{\infty}-\theta_1+\sigma)\mfS(\theta_{\infty}-\theta_1-\sigma)
 & = c_{\sigma}-c_1c_{\infty}-s_1s_{\infty} ,
 \label{trigId:c} \\
 2\mfS(\theta_{\infty}+\theta_1+\sigma)\mfS(\theta_{\infty}+\theta_1-\sigma)
 & = c_{\sigma}-c_1c_{\infty}+s_1s_{\infty} ,
 \label{trigId:d} \\
 2\mfS(\theta_{\infty}-\theta_1+\sigma)\mfS(\theta_{\infty}+\theta_1-\sigma)
 & =-c_{\infty}+c_1c_{\sigma}+s_1s_{\sigma} ,
 \label{trigId:e} \\
 2\mfS(\theta_{\infty}+\theta_1+\sigma)\mfS(\theta_{\infty}-\theta_1-\sigma)
 & =-c_{\infty}+c_1c_{\sigma}-s_1s_{\sigma} ,
 \label{trigId:f} \\
 2\mfS(\theta_0-\theta_t+\sigma)\mfS(\theta_0+\theta_t+\sigma)
 & = c_t-c_0c_{\sigma}+s_0s_{\sigma} ,
 \label{trigId:g} \\
 2\mfS(\theta_0+\theta_t-\sigma)\mfS(\theta_0-\theta_t-\sigma)
 & = c_t-c_0c_{\sigma}-s_0s_{\sigma} ,
 \label{trigId:h} \\
 2\mfS(\theta_0-\theta_t+\sigma)\mfS(\theta_0+\theta_t-\sigma)
 & =-c_0+c_tc_{\sigma}+s_ts_{\sigma} ,
 \label{trigId:i} \\
 2\mfS(\theta_0+\theta_t+\sigma)\mfS(\theta_0-\theta_t-\sigma)
 & =-c_0+c_tc_{\sigma}-s_ts_{\sigma} .
 \label{trigId:j}
\end{align}
We also note that a general member $ M $ of $ {\rm SL}(2,\CC) $, characterised by 
$ \det M=1 $ and $ {\rm tr}M =2c $, can be written as
\begin{equation}
   M = \begin{pmatrix}
       c+x & -(s+ix)b \cr
       \frac{\displaystyle s-ix}{\displaystyle b} & c-x 
       \end{pmatrix} ,
\label{SL2C_param}
\end{equation}
where $ c^2+s^2=1 $, and with arbitrary complex numbers $ b,b^{-1} \neq 0 $, $ x^{-1} \neq 0 $.
One can parameterise $ c=\cos\pi\theta, s=\sin\pi\theta $.
The observation made in \cite{Bo_2005} that simplifies the derivation is that the 
Jimbo parameterisation is constructed so that
\begin{equation}
  CM_{t}M_{0}C^{-1} 
  = \Delta := \begin{pmatrix} \epsilon_{\sigma} & 0 \cr 0 & \epsilon^{-1}_{\sigma}
                       \end{pmatrix} ,
\label{PVI_Jcon}
\end{equation}
for a non-singular matrix $ C \in {\rm GL}(2,\CC) $.
This implies
\begin{equation}
   CM^{-1}_{1}M^{-1}_{\infty}C^{-1} 
   = \begin{pmatrix} \epsilon_{\sigma} & 0 \cr 0 & \epsilon^{-1}_{\sigma}
     \end{pmatrix} .
\end{equation}
Using (\ref{PVI_Minfty}) and the parameterisation (\ref{SL2C_param}) for $ M_1 $
we can recast the above relation as a homogeneous matrix equation for $ C $, that is
\begin{equation}
     \begin{pmatrix} \epsilon_{\sigma} & 0 \cr 0 & \epsilon^{-1}_{\sigma}
     \end{pmatrix}C 
  = C\begin{pmatrix}
      \epsilon_{\infty}(c_1+x) & -\epsilon_{\infty}(s_1+ix)b \cr
      \frac{\displaystyle s_1-ix}{\displaystyle b\epsilon_{\infty}} & 
      \frac{\displaystyle c_1-x }{\displaystyle \epsilon_{\infty}}
      \end{pmatrix} .
\end{equation}
This has a solution if and only if $ ixs_{\infty}=c_{\sigma}-c_1c_{\infty} $, which fixes
$ x $. Using the identities (\ref{trigId:a},\ref{trigId:b},\ref{trigId:c},\ref{trigId:d})
the equation for the components of $ C $ can be written
\begin{align}
  -i\epsilon_{\infty}b\mfS(\theta_{\infty}+\theta_1+\sigma)C_{11}
   + \mfS(\theta_{\infty}-\theta_1-\sigma)C_{12} 
  & = 0 ,
  \label{} \\
   \epsilon_{\infty}b\mfS(\theta_{\infty}+\theta_1-\sigma)C_{21} 
   + i\mfS(\theta_{\infty}-\theta_1+\sigma)C_{22} 
  & = 0 ,
  \label{}
\end{align}
assuming $ \theta_{\infty}\pm\theta_1\pm\sigma \neq \ZZ $. If one defines an
arbitrary complex constant $ r:=ib\epsilon_{\infty} $ then this yields the solution 
(\ref{PVI_C}) for $ C $. We note that $ \det C=-s_{\infty}s_{\sigma} $. 
Employing the solution for $ x $ and the identities (\ref{trigId:c},\ref{trigId:d})
we arrive at the formula (\ref{PVI_M1}) for $ M_1 $.
For the remaining monodromy matrices we utilise the representation (\ref{SL2C_param})
for $ CM_{0}C^{-1} $ and $ CM_{t}C^{-1} $ in the relation (\ref{PVI_Jcon}). From the
$ (1,1) $ and $ (2,2) $ components of the resulting matrix equation we find
\begin{align}
  \epsilon_{\infty}c_0-c_t & = \epsilon_{\infty}x_0+x_t ,
  \label{} \\
 -\epsilon^{-1}_{\infty}c_0+c_t & = \epsilon^{-1}_{\infty}x_0+x_t ,
  \label{}
\end{align}
and one deduces that $ is_{\sigma}x_0=c_0c_{\sigma}-c_t $ and
$ is_{\sigma}x_t=c_tc_{\sigma}-c_0 $. From the other components we find that
$ x^2_t-x^2_0=c^2_t-c^2_0 $ and the formula for the ratio of the other undetermined
constants is given by
\begin{equation}
    \frac{b_t}{b_0} 
    = \epsilon_{\sigma}\frac{\mfS(\theta_0-\theta_t-\sigma)}{\mfS(\theta_0-\theta_t+\sigma)} ,
\end{equation}
where we have utilised (\ref{trigId:g},\ref{trigId:h},\ref{trigId:i},\ref{trigId:j}).
If we define the remaining undetermined constant $ b_t:=-i\epsilon_{\sigma} s_{0t} $ then
we recover the parameterisations (\ref{PVI_M0}), (\ref{PVI_Mt}) for $ M_0, M_t $.
\end{proof}

A key identity is the following connection relation which relates $ s_{0t}, \sigma_{0t} $ 
to $ \sigma_{t1} $ and $ \sigma_{01} $.
\begin{lemma}[Jimbo\cite{Ji_1982}]
One of the connection relations is
\begin{multline}
  4s_{0t}^{\pm 1}\sin\frac{\pi}{2}(\theta_{0}+\theta_{t}\mp\sigma_{0t})
            \sin\frac{\pi}{2}(\theta_{0}-\theta_{t}\pm\sigma_{0t}) \\ \times
            \sin\frac{\pi}{2}(\theta_{\infty}+\theta_{1}\mp\sigma_{0t})
            \sin\frac{\pi}{2}(\theta_{\infty}-\theta_{1}\pm\sigma_{0t}) \\
 = e^{\pm\pi i\sigma_{0t}}
   \left( \pm i\sin\pi\sigma_{0t}\cos\pi\sigma_{t1}-\cos\pi\theta_{t}\cos\pi\theta_{\infty}
          -\cos\pi\theta_{0}\cos\pi\theta_{1} \right) \\
          \pm i \sin\pi\sigma_{0t}\cos\pi\sigma_{01}+\cos\pi\theta_{t}\cos\pi\theta_{1}
          +\cos\pi\theta_{\infty}\cos\pi\theta_{0}
\label{PVI_pSoln}
\end{multline}
\end{lemma}
\begin{proof}
The proof of this has been detailed in Boalch \cite{Bo_2005} with a typographical 
correction to the original formula in \cite{Ji_1982}. So we content ourselves with
a brief summary of the steps involved. After noting (\ref{PVI_Jcon}) it was found
that the monodromy invariants are more manageable in the forms
\begin{align}
  p_{t1} & = {\rm tr}(CM^{-1}_{\infty}C^{-1}\Delta^{-1}(CM_tC^{-1})) ,
  \label{MI:a} \\
  p_{01} & = {\rm tr}(CM^{-1}_{\infty}C^{-1}\Delta^{-1}(CM_0C^{-1})) ,
  \label{MI:b}
\end{align}
because they were linear in $ 1, s_{0t}, s^{-1}_{0t} $. 
Upon taking the combination (\ref{MI:a})+$\epsilon_{\sigma}$(\ref{MI:b}) in order
to eliminate the $ s^{-1}_{0t} $ terms it
was found that the resulting expression could be factorised through the use of
the identities
\begin{multline}
  \mfS(\theta_{\infty}+\theta_1-\sigma)\mfS(\theta_{\infty}-\theta_1-\sigma) \\
 -\mfS(\theta_{\infty}-\theta_1+\sigma)\mfS(\theta_{\infty}+\theta_1+\sigma)
   = -s_{\infty}s_{\sigma} ,
\end{multline}
\begin{multline}
  \epsilon^{-1}_{\infty}
  \mfS(\theta_{\infty}+\theta_1-\sigma)\mfS(\theta_{\infty}-\theta_1-\sigma) \\
 -\epsilon_{\infty}
  \mfS(\theta_{\infty}-\theta_1+\sigma)\mfS(\theta_{\infty}+\theta_1+\sigma)
   = is_{\infty}(\epsilon_{\sigma}c_{\infty}-c_1) ,
\end{multline}
\begin{multline}
  \epsilon_{\infty}
  \mfS(\theta_{\infty}+\theta_1-\sigma)\mfS(\theta_{\infty}-\theta_1-\sigma) \\
 -\epsilon^{-1}_{\infty}
  \mfS(\theta_{\infty}-\theta_1+\sigma)\mfS(\theta_{\infty}+\theta_1+\sigma)
   = is_{\infty}(c_1-\epsilon_{\sigma}^{-1}c_{\infty}) ,
\end{multline}
enabling a factor of $ 2s_{\infty}s_{\sigma} $ to be cancelled out. This then yielded
one of the desired formulae (\ref{PVI_pSoln}), whilst the other could be found from the
other combination of (\ref{MI:a},\ref{MI:b}).
\end{proof}

A consequence of this is a constraint on the monodromy invariants 
$ \{ p_{0t},p_{t1},p_{01} \} $ which is an
algebraic variety defining a sub-manifold, the monodromy manifold, of $ \mathbb{C}^3 $.

\begin{lemma}[Jimbo\cite{Ji_1982}]
The monodromy manifold is given by
\begin{multline}
  \mfM(p_{0t},p_{t1},p_{01}) := \\
  p_{0t}p_{t1}p_{01}+p_{0t}^2+p_{t1}^2+p_{01}^2
 -(p_0p_t+p_1p_{\infty})p_{0t}-(p_tp_1+p_0p_{\infty})p_{t1}-(p_0p_1+p_tp_{\infty})p_{01}
 \\
 +p_0^2+p_t^2+p_1^2+p_{\infty}^2+p_0p_tp_1p_{\infty}-4=0 .
\label{PVI_manifold}
\end{multline}
\end{lemma}
\begin{proof}
We multiply the upper and lower sign forms of the left-hand side of (\ref{PVI_pSoln})
to eliminate $ s_{0t} $ and then employ the identities 
(\ref{trigId:a},\ref{trigId:b},\ref{trigId:g},\ref{trigId:h}) to replace the product
of sines. Equating this to the corresponding product of the right-hand sides
yields (\ref{PVI_manifold}) as the only nontrivial factor.  
\end{proof}

\begin{remark}
The above connection relation involves only the free parameters $ s_{0t} $ and 
the monodromy invariants $ \sigma_{0t},\sigma_{t1},\sigma_{01} $. As we shall see
immediately below the arbitrary parameters $ s_{0t} $ and $ \sigma_{0t} $ appear
in the expansion for the $\tau$-function  about $ t=0 $, and there exist analogous
pairs about $ t=1,\infty $. Correspondingly there exist two other forms of the
connection relation (\ref{PVI_pSoln}) involving either the parameters 
$ s_{t1}, s_{01} $ and can be deduced directly from (\ref{PVI_pSoln}) by a simple 
substitution rule given at the end of Theorem \ref{PVI_Ethm1}. However both 
connection relations yield the same formula (\ref{PVI_manifold}) for the variety
defining the monodromy manifold. 
\end{remark}

Now we come to the fundamental result for the expansion of the $\tau$-function in
the neighbourhood of the fixed singularities of the sixth Painlev\'e system
at $ t=0,1,\infty $.
\begin{theorem}[Jimbo\cite{Ji_1982}]\label{PVI_Ethm0}
Under the conditions 
(\ref{PVI_conditions:a},\ref{PVI_conditions:b},\ref{PVI_conditions:c})
we have the expansion of the $\tau$-function as $ t \to 0 $ in the domain
$ \{t\in \CC| 0<|t|<\varepsilon, |{\rm arg}(t)|<\phi\} $ for all $ \varepsilon>0 $ 
and any $ \phi>0 $
\begin{multline}
  \tau(t) \sim Ct^{(\sigma^2-\theta^2_0-\theta^2_t)/4} \\
   \times\Bigg\{ 1 
               + \frac{(\theta^2_0-\theta^2_t-\sigma^2)(\theta^2_{\infty}-\theta^2_1-\sigma^2)}
                     {8\sigma^2}t \\
   - \hat{s}\frac{[\theta^2_0-(\theta_t-\sigma)^2][\theta^2_{\infty}-(\theta_1-\sigma)^2]}
                 {16\sigma^2(1+\sigma)^2}t^{1+\sigma} \\
   - \hat{s}^{-1}\frac{[\theta^2_0-(\theta_t+\sigma)^2][\theta^2_{\infty}-(\theta_1+\sigma)^2]}
                 {16\sigma^2(1-\sigma)^2}t^{1-\sigma}
   + {\rm O}(|t|^{2(1-\Re(\sigma))}) \Bigg\}
\label{PVI_tExp_0}
\end{multline}
where $ \sigma \neq 0 $ and $ \hat{s} $ are related to $ s $ through
\begin{multline}
  \hat{s} = s \frac{\Gamma^2(1-\sigma)\Gamma(1+\frac{1}{2}(\theta_0+\theta_t+\sigma))
                    \Gamma(1+\frac{1}{2}(-\theta_0+\theta_t+\sigma))}
                   {\Gamma^2(1+\sigma)\Gamma(1+\frac{1}{2}(\theta_0+\theta_t-\sigma))
                    \Gamma(1+\frac{1}{2}(-\theta_0+\theta_t-\sigma))} \\
        \times\frac{\Gamma(1+\frac{1}{2}(\theta_{\infty}+\theta_1+\sigma))
                    \Gamma(1+\frac{1}{2}(-\theta_{\infty}+\theta_1+\sigma))}
                   {\Gamma(1+\frac{1}{2}(\theta_{\infty}+\theta_1-\sigma))
                    \Gamma(1+\frac{1}{2}(-\theta_{\infty}+\theta_1-\sigma))} ,
\label{PVI_s_0t}
\end{multline}
and we employ the short-hand notation $ s=s_{0t} $, $ \hat{s}=\hat{s}_{0t} $ and 
$ \sigma=\sigma_{0t} $. The monodromy data defining the unique solution to the sixth
Painlev\'e system is $ \{\sigma_{0t}, s_{0t}\} $.
Here $ C $ is an arbitrary constant.
\end{theorem}

\begin{proof}
The details of this proof of this are given in Jimbo \cite{Ji_1982} and so we do not
repeat them here. Also Guzzetti has laid out some of the intermediate steps in the
appendix of his work on the elliptic representations of the general Painlev\'e
six equation \cite{Gu_2002}. 
\end{proof}

The regular singularities $ x=0,t,1,\infty $ play equivalent roles and can be exchanged
under linear fractional or M\"obius transformations. Consequently one can solve the
connection problem very neatly and under the additional conditions
\begin{align}
  & 0 < \Re(\sigma_{t1}),\Re(\sigma_{01}) < 1 ,
  \label{PVI_conditions:d}\\
  & \theta_1\pm\theta_t\pm\sigma_{t1},\quad 
    \theta_{\infty}\pm\theta_0\pm\sigma_{t1} \notin 2\mathbb{Z} ,
  \label{PVI_conditions:e}\\
  & \theta_0\pm\theta_1\pm\sigma_{01},\quad 
    \theta_{\infty}\pm\theta_t\pm\sigma_{01} \notin 2\mathbb{Z} ,
  \label{PVI_conditions:f}
\end{align}
derive expansions about $ t=1,\infty $.

\begin{theorem}[Jimbo\cite{Ji_1982}]\label{PVI_Ethm1}
Under the conditions 
(\ref{PVI_conditions:a},\ref{PVI_conditions:d},\ref{PVI_conditions:e})
we have the expansion of the $\tau$-function as $ t \to 1 $
\begin{multline}
  \tau(t) \sim C(1-t)^{(\sigma_{t1}^2-\theta^2_1-\theta^2_t)/4} \\
   \times\Bigg\{ 1 
               + \frac{(\theta^2_1-\theta^2_t-\sigma_{t1}^2)
                       (\theta^2_{\infty}-\theta^2_0-\sigma_{t1}^2)}
                     {8\sigma_{t1}^2}(1-t) \\
   - \hat{s}_{t1}\frac{[\theta^2_1-(\theta_t-\sigma_{t1})^2]
                  [\theta^2_{\infty}-(\theta_0-\sigma_{t1})^2]}
                 {16\sigma_{t1}^2(1+\sigma_{t1})^2}(1-t)^{1+\sigma_{t1}} \\
   - \hat{s}_{t1}^{-1}\frac{[\theta^2_1-(\theta_t+\sigma_{t1})^2]
                       [\theta^2_{\infty}-(\theta_0+\sigma_{t1})^2]}
                 {16\sigma_{t1}^2(1-\sigma_{t1})^2}(1-t)^{1-\sigma_{t1}}
   + {\rm O}(|1-t|^{2(1-\Re(\sigma_{t1}))}) \Bigg\} ,
\label{PVI_tExp_1}
\end{multline}
and as $ t \to \infty $
\begin{multline}
  \tau(t) \sim Ct^{-(\sigma_{01}^2-\theta^2_{\infty}+\theta^2_t)/4} \\
   \times\Bigg\{ 1 
               + \frac{(\theta^2_{\infty}-\theta^2_t-\sigma_{01}^2)
                       (\theta^2_{0}-\theta^2_1-\sigma_{01}^2)}
                     {8\sigma_{01}^2}t^{-1} \\
   - \hat{s}_{01}\frac{[\theta^2_{\infty}-(\theta_t-\sigma_{01})^2]
                       [\theta^2_{0}-(\theta_1-\sigma_{01})^2]}
                 {16\sigma_{01}^2(1+\sigma_{01})^2}t^{-1-\sigma_{01}} \\
   - \hat{s}_{01}^{-1}\frac{[\theta^2_{\infty}-(\theta_t+\sigma_{01})^2]
                            [\theta^2_{0}-(\theta_1+\sigma_{01})^2]}
                 {16\sigma_{01}^2(1-\sigma_{01})^2}t^{-1+\sigma_{01}}
   + {\rm O}(|t|^{-2(1-\Re(\sigma_{01}))}) \Bigg\} .
\label{PVI_tExp_infty}
\end{multline}
Here $ \hat{s}_{t1}, \hat{s}_{01} $ are found by making the following substitutions in
(\ref{PVI_s_0t}), (\ref{PVI_pSoln}) respectively
\begin{align}
   & \hat{s} \to \hat{s}_{t1}, \quad s \to s_{t1}, \quad \theta_{0} \leftrightarrow \theta_{1},
     \quad \sigma \to \sigma_{t1}, \quad \sigma_{t1} \to \sigma_{0t} ,
   \\
   & \hat{s} \to \hat{s}_{01}, \quad s \to s_{01}, 
     \quad \theta_{0} \leftrightarrow \theta_{\infty},
     \quad \sigma \to \sigma_{01}, \quad \sigma_{01} \to \tilde{\sigma}_{01} ,
\end{align}
with
\begin{equation}
  \cos\pi\tilde{\sigma}_{01} = -\cos\pi\sigma_{0t}-2\cos\pi\sigma_{01}\cos\pi\sigma_{t1}
                        +2(\cos\pi\theta_0\cos\pi\theta_t+\cos\pi\theta_{\infty}\cos\pi\theta_1) .
\end{equation}
The monodromy data defining the unique solution to the sixth Painlev\'e system is 
either $ \{\sigma_{t1}, s_{t1}\} $ or $ \{\sigma_{01}, s_{01}\} $.
\end{theorem}

\section{Monodromy Data for the Spectrum Singularity Ensemble}\label{sectionMData_SSE}
\setcounter{equation}{0}

The precise relationship between the spectrum singularity average $ A_N(t;) $ and the 
isomonodromy theory of the sixth Painlev\'e system is given by the following result.
Its validity relies on the conjecture that the expansions of Jimbo
given in Theorems \ref{PVI_Ethm0} and \ref{PVI_Ethm1} remain valid upon relaxation of the
constraints (\ref{PVI_conditions:a},\ref{PVI_conditions:b},\ref{PVI_conditions:c}) and 
(\ref{PVI_conditions:d},\ref{PVI_conditions:e},\ref{PVI_conditions:f}), 
provided the former are well defined (i.e. do not then diverge).

\begin{proposition}
For the spectrum singularity ensemble the associated isomonodromic system is not unique
but the monodromy data for any of these systems falls into three generic cases.
An example of each case is given below in cases (A), (B) and (C). The formal monodromy 
exponents can be taken to belong to either of three sets
\begin{align}
   {\rm Case}(A): &\quad
   \theta_{0} = -\mu-\omega ,\quad 
   \theta_{t} = N+2\omega_1 ,\quad
   \theta_{1} = N+2\mu ,\quad
   \theta_{\infty} =-\mu-\bar{\omega} ,
   \label{SSE_Mexp:A} \\
   {\rm Case}(B): &\quad
   \theta_{0} = \mu-\bar{\omega} ,\quad 
   \theta_{t} = N ,\quad
   \theta_{1} = N+2\mu+2\omega_1 ,\quad
   \theta_{\infty} = \mu-\omega , 
   \label{SSE_Mexp:B} \\
   {\rm Case}(C): &\quad
   \theta_{0} = -2\omega_1 ,\quad 
   \theta_{t} = N+\mu+\omega ,\quad
   \theta_{1} = N+\mu+\bar{\omega} ,\quad
   \theta_{\infty} = 2\mu . 
   \label{SSE_Mexp:C}
\end{align}
The monodromy invariants for either case are
\begin{equation}
   \sigma_{0t}=N-\mu+\bar{\omega}, \quad
   \sigma_{t1}=2\mu+2\omega_1, \quad
   \sigma_{01}=N-\mu+\omega .
\label{SSE_MI}
\end{equation}
In the case (A) the monodromy coefficients are
\begin{align}
  s_{0t}
  & = 1+\frac{2i\sin\pi(\mu-\bar{\omega})}{\xi^*e^{-\pi i(\mu-\bar{\omega})}} ,
  \label{SSE_s0t:A} \\
  s_{t1}
  & = 1+\xi^*\frac{e^{-\pi i(\mu-\bar{\omega})}}{2i}
      \frac{\sin\pi(2\mu+2\omega_1)}{\sin\pi 2\mu\sin\pi(\mu+\omega)} ,
  \label{SSE_st1:A}  \\
  s_{01}
  & = -\frac{\xi^*-1+e^{2\pi i(\mu+\omega)}}
            {\xi^*-1+e^{4\pi i\mu}} .
  \label{SSE_s01:A} 
\end{align}
All monodromy matrices are lower triangular
\begin{gather}
  M_0 = \begin{pmatrix}
          e^{-\pi i(\mu+\omega)} & 0 \\ m_0 & e^{\pi i(\mu+\omega)} \\
        \end{pmatrix} ,
  \label{SSE_MA:a}\\
  M_t = \begin{pmatrix}
          e^{\pi i(N+2\omega_1)} & 0 \\ m_t  & e^{-\pi i(N+2\omega_1)} \\
        \end{pmatrix} ,
  \label{SSE_MA:b}\\
  M_1 = \begin{pmatrix}
          e^{\pi i(N+2\mu)} & 0 \\ m_1 & e^{-\pi i(N+2\mu)} \\
        \end{pmatrix} ,
  \label{SSE_MA:c}
\end{gather}
where
\begin{align}
  m_0 &=\frac{2i}{\sin\pi(\mu-\bar{\omega})}
  \left\{  \frac{\sin\pi2\omega_1\sin\pi(\mu+\bar{\omega})}{s_{0t}}
          -\frac{\sin\pi2\mu\sin\pi(\mu+\omega)}{r} \right\} ,
  \label{SSE_MA:d}\\
  m_t &=\frac{2i(-1)^N\sin\pi2\omega_1}{\sin\pi(\mu-\bar{\omega})}
  \left\{ -\frac{\sin\pi(\mu+\bar{\omega})}{s_{0t}}e^{\pi i(\mu-\bar{\omega})}
          +\frac{\sin\pi 2\mu}{r} \right\} ,
  \label{SSE_MA:e}\\
  m_1 &=-\frac{2i(-1)^N\sin\pi 2\mu}{r}e^{-\pi i(\mu+\bar{\omega})} .
  \label{SSE_MA:f}
\end{align}

\begin{sideways}
\begin{minipage}[c][6cm][c]{17cm}{
For Case (B) the monodromy coefficients are
\begin{align}
  s_{0t}
  & = 1+\frac{2i\sin\pi(\mu-\bar{\omega})}{\xi^*e^{-\pi i(\mu-\bar{\omega})}} ,
  \label{SSE_s0t:B} \\
  s_{t1}
  \frac{\sin\pi 2\omega_1\sin\pi(\mu+\bar{\omega})}{\sin\pi(2\mu+2\omega_1)}
  & = \frac{\sin\pi 2\mu\sin\pi(\mu+\omega)}{\sin\pi(2\mu+2\omega_1)}
      +\xi^*\frac{e^{-\pi i(\mu-\bar{\omega})}}{2i} ,
  \label{SSE_st1:B}  \\
  s_{01}
  & = -\frac{\xi^*-1+e^{2\pi i(\mu+\omega)}}
            {\xi^*-1+e^{4\pi i\mu}} .
  \label{SSE_s01:B} 
\end{align}
One of monodromy matrices is proportional to the identity, the others are full
\begin{gather}
  M_0 = \frac{i}{\sin\pi(\mu-\omega)}
        \begin{pmatrix}
          e^{-\pi i(\mu-\omega)}\cos\pi(\mu-\bar{\omega})-\cos\pi(2\mu+2\omega_1) &
          2r\sin\pi(\mu+\bar{\omega})\sin\pi2\mu  \\
          -\frac{\displaystyle 2}{\displaystyle r}\sin\pi(\mu+\omega)\sin\pi2\omega_1 &
          -e^{\pi i(\mu-\omega)}\cos\pi(\mu-\bar{\omega})+\cos\pi(2\mu+2\omega_1) \\
        \end{pmatrix} ,
  \label{SSE_MB:a}\\
  M_t = (-1)^N I ,
  \label{SSE_MB:b}\\
  M_1 = \frac{i(-1)^N}{\sin\pi(\mu-\omega)}
        \begin{pmatrix}
          e^{-\pi i(\mu-\omega)}\cos\pi(2\mu+2\omega_1)-\cos\pi(\mu-\bar{\omega}) &
          2re^{-\pi i(\mu-\omega)}\sin\pi(\mu+\bar{\omega})\sin\pi2\mu  \\
          \frac{\displaystyle 2}{\displaystyle r}
               e^{\pi i(\mu-\omega)}\sin\pi(\mu+\omega)\sin\pi2\omega_1 &
          \cos\pi(\mu-\bar{\omega})-e^{\pi i(\mu-\omega)}\cos\pi(2\mu+2\omega_1) \\
        \end{pmatrix} .
  \label{SSE_MB:c}
\end{gather}}
\end{minipage}
\end{sideways}

For Case (C) the monodromy coefficients are
\begin{align}
  s_{0t}
  & = 1+\frac{2i\sin\pi(\mu-\bar{\omega})}{\xi^*e^{-\pi i(\mu-\bar{\omega})}} ,
  \label{SSE_s0t:C} \\
  s_{t1}
  \frac{\sin\pi 2\omega_1\sin\pi(\mu+\bar{\omega})}{\sin\pi(2\mu+2\omega_1)}
  & = \frac{\sin\pi 2\mu\sin\pi(\mu+\omega)}{\sin\pi(2\mu+2\omega_1)}
      +\xi^*\frac{e^{-\pi i(\mu-\bar{\omega})}}{2i} ,
  \label{SSE_st1:C}  \\
  s_{01}
  & = -\frac{\xi^*-1+e^{2\pi i(\mu+\omega)}}
            {\xi^*-1+e^{4\pi i\mu}} .
  \label{SSE_s01:C} 
\end{align}
All monodromy matrices are upper triangular
\begin{gather}
  M_0 = \begin{pmatrix}
          e^{\pi i2\omega_1} & m_0 \\ 0 & e^{-\pi i2\omega_1} \\
        \end{pmatrix} ,
  \label{SSE_MC:a}\\
  M_t = \begin{pmatrix}
          e^{-\pi i(N+\mu+\omega)} & m_t \\ 0 & e^{\pi i(N+\mu+\omega)} \\
        \end{pmatrix} ,
  \label{SSE_MC:b}\\
  M_1 = \begin{pmatrix}
          e^{-\pi i(N+\mu+\bar{\omega})} & m_1 \\ 0 & e^{\pi i(N+\mu+\bar{\omega})} \\
        \end{pmatrix} ,
  \label{SSE_MC:c}
\end{gather}
where
\begin{align}
  m_0 &=\frac{2i}{\sin\pi(\mu-\bar{\omega})}
  \left\{ -\sin\pi2\mu\sin\pi(\mu+\omega)s_{0t}
          +\sin\pi2\omega_1\sin\pi(\mu+\bar{\omega})r \right\} ,
  \label{SSE_MC:d}\\
  m_t &=\frac{2i(-1)^N\sin\pi(\mu+\omega)}{\sin\pi(\mu-\bar{\omega})}
  \left\{  \sin\pi 2\mu e^{-\pi i(\mu-\bar{\omega})}s_{0t}
          -\sin\pi(\mu+\bar{\omega})r \right\} ,
  \label{SSE_MC:e}\\
  m_1 &= 2i\sin\pi(\mu+\bar{\omega}) e^{-\pi i(N+2\mu)}r .
  \label{SSE_MC:f}
\end{align}
\end{proposition}

\begin{proof}
Comparison of the two differential equations for the $\sigma$-function, 
(\ref{PVI_sigmaF}) and (\ref{PVI_SF}), imply that in general
\begin{equation} 
  \{v_1,v_2,v_3,v_4\} = \frac{1}{2}
  \begin{cases}
   \epsilon_1 (\theta_{t}+\theta_{\infty}) \\
   \epsilon_2 (\theta_{t}-\theta_{\infty}) \\
   \epsilon_3 (\theta_{0}+\theta_{1}) \\
   \epsilon_4 (\theta_{0}-\theta_{1}) \\
  \end{cases} ,
\end{equation} 
with $ \epsilon_j=\pm 1, j=1,2,3,4 $ and $ \epsilon_1\epsilon_2\epsilon_3\epsilon_4=1 $. 
Using either set of parameters, (\ref{AN_paramJUE}) from the JUE correspondence or 
(\ref{AN_paramCyUE}) from the CyUE correspondence, we find that the monodromy 
exponents can be given by one of three sets
\begin{equation}
  \{ \theta_{0},\theta_{t},\theta_{1},\theta_{\infty}\}
 =   \begin{cases}
        N+2\mu, N+2\omega_1, \mu+\omega, \mu+\bar{\omega} \\
        N, N+2\mu+2\omega_1, \mu-\omega, \mu-\bar{\omega} \\
        N+\mu+\omega, N+\mu+\bar{\omega}, -2\mu, 2\omega_1
     \end{cases} ,
\end{equation}
modulo permutations of the monodromy exponents and an even number of sign reversals. 
This is a manifestation of the non-uniqueness of the isomonodromic system for our
problem.
For definiteness we choose one example of the three cases, namely the cases
(A), (B) and (C) given in 
(\ref{SSE_Mexp:A},\ref{SSE_Mexp:B},\ref{SSE_Mexp:C}).
We note some simple identities which do not depend on the choice of the permutation 
or the sign
\begin{gather}
   v_1^2+v_2^2+v_3^2+v_4^2 =  \frac{1}{2}(\theta_0^2+\theta_t^2+\theta_1^2+\theta_{\infty}^2) ,
   \label{Mexp_id:a}\\
   v_1v_2v_3v_4 = \frac{1}{16}(\theta_0^2-\theta_1^2)(\theta_t^2-\theta_{\infty}^2) . 
   \label{Mexp_id:b}
\end{gather}
and in particular the following products which apply equally to cases (A), (B) and (C)
\begin{gather}
  (\theta_0^2-\theta_1^2)(\theta_{\infty}^2-\theta_t^2) =
  (N-\mu+\omega)(N+\mu-\omega)(N+3\mu+\omega)(N+\mu+2\bar{\omega}+\omega) ,
  \\
  (\theta_0^2-\theta_t^2)(\theta_{\infty}^2-\theta_1^2) =
  (N-\mu+\bar{\omega})(N+\mu-\bar{\omega})(N+3\mu+\bar{\omega})(N+\mu+2\omega+\bar{\omega}) ,
  \\
  (\theta_1^2-\theta_t^2)(\theta_{\infty}^2-\theta_0^2) =
  (2N+2\mu+2\omega_1)(2\mu-2\omega_1)(2\mu+2\omega_1)(\bar{\omega}-\omega) .
\end{gather}
If we make a comparison of the $\tau$-functions themselves for the JUE correspondence, 
(\ref{PVI_sigma}) and (\ref{SSE_PVIsoln}) we find, at the level of the 
$\tau$-functions
\begin{multline}
   A_N(t;) = \\
   \tilde{C}t^{\frac{1}{8}(\theta_0^2+\theta_t^2-\theta_1^2-\theta_{\infty}^2)
                        -\frac{1}{2}e_2[v^{\rm JUE}]-\mu N}
             (1-t)^{\frac{1}{8}(-\theta_0^2+\theta_t^2+\theta_1^2-\theta_{\infty}^2)
                        -e'_2[v^{\rm JUE}]+\frac{1}{2}e_2[v^{\rm JUE}]}
   \tau(t;\theta).
\label{AN_tauFn}
\end{multline}

Now applying the $ t\to 0 $ expansions for $ \tau(t;\theta) $, namely 
(\ref{PVI_tExp_0}),
with those of $ A_N(t;) $, (\ref{AN_exp_0}), we first note that the exponent of the 
$t$-prefactors must be consistent and this implies
\begin{equation}
   \sigma_{0t}^2 = \frac{1}{2}(\theta_0^2+\theta_t^2+\theta_1^2+\theta_{\infty}^2)
                   +2e_2[v^{\rm JUE}]
   = (\sum_{i=1}^{4}v^{\rm JUE}_i)^2 = (N-\mu+\bar{\omega})^2 .
\end{equation}
This result applies to all the cases and for definiteness we make the choice of sign
outlined in (\ref{SSE_MI}) (the other choice of sign is essentially equivalent).
Turning to the leading analytic term of order $ t $ in both expansions we find that
its coefficient is given precisely by
\begin{equation}
   \frac{2\mu(\mu+\omega)N}{N-\mu+\bar{\omega}} ,
\end{equation} 
for all three cases upon employing our solution for $ \sigma_{0t}^2 $.
Next we make a comparison of the non-analytic terms in (\ref{AN_exp_0}) and 
(\ref{PVI_tExp_0}) and it is here that we have to treat the cases separately.
However it is generally true that
in the classical situation of the finite rank random matrix ensemble only one
of the non-analytic terms is ever present, the other being switched off through
the following mechanism. 
Taking case (A) first our task is to show how the coefficient of the non-analytic 
term $ t^{1+\sigma_{0t}} $ vanishes and to match the remaining coefficient with that 
of the application. The identification of $ \sigma_{0t} $ immediately implies 
$ \sigma_{0t}=\theta_{0}+\theta_{t} $. Therefore
\begin{equation}
  [t^{1+\sigma_{0t}}] \propto \frac{1}{\Gamma(\frac{\theta_{0}+\theta_{t}-\sigma_{0t}}{2})}
   = 0 .
\end{equation}  
From a comparison of the coefficients of the remaining non-analytic term we find 
precise agreement and this enables us to determine the solution 
for the monodromy coefficient $ s_{0t} $ in (\ref{SSE_s0t:A}).
For case (B) we note that $ \sigma_{0t}=\theta_{t}-\theta_{0} $ and thus
\begin{equation}
  [t^{1+\sigma_{0t}}] \propto \frac{1}{\Gamma(\frac{-\theta_{0}+\theta_{t}-\sigma_{0t}}{2})}
   = 0 .
\end{equation}
Similarly we find agreement for the coefficient of the surviving non-analytic term 
and that enables us to fix the monodromy coefficient as in (\ref{SSE_s0t:B}).
In case (C) we have  $ \sigma_{0t}=\theta_{1}-\theta_{\infty} $ and thus
\begin{equation}
  [t^{1+\sigma_{0t}}] \propto \frac{1}{\Gamma(\frac{-\theta_{\infty}+\theta_{1}-\sigma_{0t}}{2})}
   = 0 .
\end{equation}
Again we find agreement for the coefficient of the $ t^{1+\sigma_{0t}} $ term
and conclude that the monodromy coefficient is given by (\ref{SSE_s0t:C}).
For each of the three cases we observe that the Jimbo parameterisation of the 
monodromy fails 
(see condition (\ref{PVI_conditions:c}) however a meaningful result emerges so we
conjecture that the Theorem still holds with these relaxations.
 
Now we make a comparison of the expansions at $ t=1 $, namely
(\ref{AN_exp_1}) and (\ref{PVI_tExp_1}). Examination of the algebraic prefactor 
using the relation (\ref{AN_tauFn}) leads us to conclude 
$ \sigma_{t1}^2 = (2\mu+2\omega_1)^2 $
in all three cases and we choose the positive sign, as in (\ref{SSE_MI}).
Employing this solution we compute the coefficient of the analytic term is, in 
all three cases,
\begin{equation}
  \frac{\mu(\bar{\omega}-\omega)N}{2\mu+2\omega_1} ,
\end{equation}
which is are entirely consistent with that in (\ref{AN_exp_1}). To examine the 
non-analytic terms we take the three cases separately again. For case (A) we see 
that the coefficient of the $ (1-t)^{1-\sigma_{t1}} $ vanishes because 
$ \sigma_{t1}=-\theta_0-\theta_{\infty} $ and 
\begin{equation}
  [(1-t)^{1-\sigma_{t1}}] \propto \frac{1}{\Gamma(\frac{\theta_{\infty}+\theta_{0}+\sigma_{t1}}{2})}
   = 0 .
\end{equation}
The coefficients of $ (1-t)^{1+\sigma_{t1}} $ now agree precisely provided we 
have the solution (\ref{SSE_st1:A}) for the monodromy coefficient $ s_{t1} $.
For case (B) we have the relation $ \sigma_{t1}=\theta_1-\theta_{t} $ and see that
\begin{equation}
  [(1-t)^{1-\sigma_{t1}}] \propto \frac{1}{\Gamma(\frac{-\theta_{1}+\theta_{t}+\sigma_{t1}}{2})}
   = 0 .
\end{equation}
Again the coefficients of $ (1-t)^{1+\sigma_{t1}} $ agree and the solution 
(\ref{SSE_st1:B}) for $ s_{t1} $ follows. In case (C) we see that 
$ \sigma_{t1}=\theta_{\infty}-\theta_{0} $ and this ensures
\begin{equation}
  [(1-t)^{1-\sigma_{t1}}] \propto \frac{1}{\Gamma(\frac{-\theta_{\infty}+\theta_0+\sigma_{t1}}{2})}
   = 0 .
\end{equation}
In this case we also find the coefficients of the remaining non-analytic terms are 
precisely consistent, leading us to deduce the solution (\ref{SSE_st1:C}) for 
$ s_{t1} $.

It remains to make a 
comparison of the expansions at $ t=\infty $, namely (\ref{AN_exp_infty}) and 
(\ref{PVI_tExp_infty}). Using (\ref{AN_tauFn}) we see the 
algebraic prefactor implies that $ \sigma_{01}^2 = (N-\mu+\omega)^2 $, and we
choose the positive sign for the exponent. Using this value for $ \sigma_{01} $ we 
compute that the coefficient for the $ t^{-1} $ term in all three cases is
\begin{equation}
  \frac{2\mu(\mu+\bar{\omega})N}{N-\mu+\omega} ,
\end{equation}
which is consistent with (\ref{AN_exp_infty}).
To treat the non-analytic terms we take the cases separately. For case (A) we note 
that $ \sigma_{01}=\theta_{t}+\theta_{\infty} $ and this implies
\begin{equation}
  [t^{-1-\sigma_{01}}] \propto \frac{1}{\Gamma(\frac{\theta_{\infty}+\theta_t-\sigma_{01}}{2})}
   = 0 .
\end{equation}
The coefficient of the $ t^{-1+\sigma_{01}} $ term is found to be in agreement with
that of (\ref{AN_exp_infty}) if we take the solution (\ref{SSE_s01:A}) for $ s_{01} $.
For case (B) the relation is $ \sigma_{01}=\theta_{t}-\theta_{\infty} $ and this
in turn implies
\begin{equation}
  [t^{-1-\sigma_{01}}] \propto \frac{1}{\Gamma(\frac{-\theta_{\infty}+\theta_t-\sigma_{01}}{2})}
   = 0 .
\end{equation} 
Again exact agreement is found for the other coefficient provided that 
(\ref{SSE_s01:B})
holds. Lastly in case (C) we have the same relation as above and the absence of the
$ t^{-1-\sigma_{01}} $ term. Examination of the coefficients of $ t^{-1+\sigma_{01}} $ 
then lead us to the solution (\ref{SSE_s01:C}).

Now we come to consideration of the connection relation (\ref{PVI_pSoln}) for 
$ t=0 $ and its two equivalent forms for $ t=1,\infty $ with respect to our solutions 
for the monodromy data.
We compute that the three connection relations of either sign decouple into a 
left-hand side and a right-hand side which vanish separately for all the cases (A), 
(B) and (C). The
left-hand sides for the $ t=0 $ connection relation vanish because
$ \theta_0+\theta_t-\sigma_{0t}=0 $ and $ \theta_{\infty}+\theta_1+\sigma_{0t}=2N $
for case (A),
$ \theta_0-\theta_t+\sigma_{0t}=0 $ and $ \theta_{0}+\theta_t+\sigma_{0t}=2N $
for case (B), and
$ \theta_{\infty}-\theta_1+\sigma_{0t}=0 $ and $ \theta_{0}+\theta_t+\sigma_{0t}=2N $
for case (C). 
Similar reasoning applies to the connection relations at $ t=1 $ and $ t=\infty $.
The right-hand sides of the relations vanish identically for both signs with the 
evaluations
of $ \sigma_{0t},\sigma_{t1},\sigma_{01} $ as given in (\ref{SSE_MI}).

To conclude we compute the monodromy matrices for the three cases and note that 
case (A), case (B) and case (C) yield the classical monodromy structure of
lower triangular matrices, full matrices with one being a signed multiple of the
identity, and upper triangular matrices respectively.
\end{proof}

\begin{remark}
Our results are consistent with the findings of Mazzocco \cite{Mz_2002} which state
that the classical non-algebraic solutions for \PVI have either reducible monodromy
groups (cases (A) and (C)) or at least one monodromy matrix is equal to 
$ \pm I $, that is the monodromy group is $1$-smaller (case (B)). Both these cases
cover the situation of a one parameter ($ N $) family of classical solutions.
\end{remark}

\begin{remark}
We observe that the exponents $ \sigma_{0t},\sigma_{t1},\sigma_{01} $ are not free 
boundary conditions for classical solutions but are fixed by certain combinations of
the formal monodromy exponents. Related to this phenomenon is that all the connection
relations decouple so that the coefficients of the monodromy coefficients
$ s_{0t},s_{t1},s_{01} $ all vanish and thus cannot be determined 
from these relations. There is a geometrical picture of the classical solutions,
which was discussed in relation to Painlev\'e II by Its and Kapaev \cite{IK_2000}.
The classical solutions of \PVI define singular points in the monodromy manifold
which are characterised by $ \mfM= 0 $ and 
\begin{align}
 \frac{\partial}{\partial p_{0t}}\mfM 
      & = p_{t1}p_{01}+2p_{0t}-p_{0}p_{t}-p_{1}p_{\infty} = 0, \\
 \frac{\partial}{\partial p_{t1}}\mfM 
      & = p_{0t}p_{01}+2p_{t1}-p_{t}p_{1}-p_{0}p_{\infty} = 0, \\
 \frac{\partial}{\partial p_{01}}\mfM 
      & = p_{0t}p_{t1}+2p_{01}-p_{0}p_{1}-p_{t}p_{\infty} = 0 .
\end{align}
We verify that these relations are satisfied for the cases (A), (B) and (C).
\end{remark}

\section*{Acknowledgements}
This work was supported by the Australian Research Council.

\bibliographystyle{plain}
\bibliography{moment,random_matrices,nonlinear}

\end{document}